\documentclass[leqno]{article}%
\usepackage{amsfonts}
\usepackage{amsmath}
\usepackage{amssymb}
\usepackage{graphicx}%
\providecommand{\U}[1]{\protect\rule{.1in}{.1in}}
\newtheorem{theorem}{Theorem}

\newtheorem{definition}[theorem]{Definition}

\newtheorem{lemma}[theorem]{Lemma}

\newtheorem{proposition}[theorem]{Proposition}
\newtheorem{remark}[theorem]{Remark}

\newenvironment{proof}[1][Proof]{\noindent\textbf{#1.} }{\ \rule{0.5em}{0.5em}}
\begin{document}

\title{A Generalized Mixed Zero-sum Stochastic Differential Game and Double Barrier
Reflected BSDEs with Quadratic Growth Coefficient}
\author{Sa\"{\i}d Hamad\`{e}ne$^{\text{ }a}$, Eduard Rotenstein$^{\text{ }b,}%
$\thanks{The work of this author was done during a visit to the Universit\'{e}
du Maine, whose generous support is gratefully acknowledged; the visit was
supported by the Grant ID\_395/2007 and CNCSIS 1156/2006.} , Adrian
Z\u{a}linescu$^{\text{ }c,\ast}$\bigskip\\$^{a\text{ }}${\small Laboratoire de Statistique et Processus, Universit\'{e}
du Maine, 72085, Le Mans Cedex 9,}\\{\small France, e-mail: hamadene@univ-lemans.fr}\\$^{b\text{ }}${\small Department of Mathematics, "Al. I. Cuza" University, Bd.
Carol no. 9-11, Ia\c{s}i,} \\{\small Rom\^{a}nia, e-mail: eduard.rotenstein@uaic.ro}\\$^{c\text{ }}${\small Institut of Mathematics of the Romanian Academy,
Ia\c{s}i branch, Bd. Carol no. 8,}\\{\small Rom\^{a}nia, e-mail: adrian.zalinescu@googlemail.com}}
\date{}
\maketitle

\begin{abstract}
This article is dedicated to the study of mixed zero-sum two-player stochastic
differential games in the situation when the player's cost functionals are
modeled by doubly controlled reflected backward stochastic equations with two
barriers whose coefficients have quadratic growth in $Z$. This is a
generalization of the risk-sensitive payoffs. We show that the lower and the
upper value function associated with this stochastic differential game with
reflection are deterministic and they are also the unique viscosity solutions
for two Isaacs equations with obstacles. \medskip

\end{abstract}

\medskip

\textbf{AMS Classification subjects: }60G40, 60H10, 60H99, 90C39, 93E05, 91A15.

\section{Introduction}

An American option is a contract which gives the right to its holder to
exercise, i.e. to ask for the wealth, when he decides before the maturity of
the option. The liabilities of the seller of the option is to provide this
wealth and nothing else. In recent last years there have been other options,
which have attracted a lot of research activity, which look like to American
options but they give also the right to the seller of the option to recall it
if he accepts to pay a money penalty, namely American game options. Those
types of options, introduced by Y.Kifer in \cite{yk} who terms them as Israeli
options, have been well documented in several papers of which we can quote
\cite{kk}, \cite{bk}, \cite{H}, \cite{kk2}. The main reason for the
introduction of those options is that brokers face high level risky
environment/market and therefore introduce clauses which allow them to
withdraw from their liabilities in case when for example a stock or a
commodity, such as crude oil, copper, steel,..., becomes more and more
expensive. In order to tackle American options we mainly use Snell envelope of
processes while in order to do so for the American game options we use value
functions of zero-sum Dynkin games. In the standard Black and Scholes model,
the value of the American game option is given by the value of the zero-sum
Dynkin game under the risk neutral probability. Additionally a hedging
strategy for the seller of the option exists (see e.g. \cite{yk,H,kk} for more details).

Now let us introduce a specific stochastic Ramsey model in a growth model in
finance (see e.g. \cite{ab} for more details). Assume we have a capital from
which is withdrawn a consumption and whose dynamics is given by:
\[
dX_{t}=X_{t}(r-c_{t})dt+X_{t}\sigma dW_{t},\,t\leq T\mbox{ and
}X_{0}=x>0.
\]
Here $r$ is the spot mean-return of the capital and $c$ is a proportion of the
capital which is consumed. Therefore usually the main objective is to find an
optimal consumption process with respect to an index which indicates the
satisfaction of the capital holder. This index depends of course on
$(c_{t})_{t\leq T}$ but also on $X$ and many other parameters such as risk
sensitiveness, utility and so on. \medskip

The problem we consider in this paper is in a way a generalization of the
combination of the two previous ones. Additionally we assume that the
criterion is of risk-sensitive type. Actually assume we have a system on which
intervene two agents $a_{1}$ and $a_{2}$ and whose dynamics is given by:
\begin{equation}
\left\{
\begin{array}
[c]{l}%
dX_{s}^{\alpha,\beta}=b(s,X_{s}^{\alpha,\beta},\alpha_{s},\beta_{s}%
)ds+\sigma(s,X_{s}^{\alpha,\beta},\alpha_{s},\beta_{s})dW_{s},\ s\in\left[
0,T\right]  ;\smallskip\\
X_{0}^{\alpha,\beta}=x.
\end{array}
\right.  \label{SDEintro}%
\end{equation}
The stochastic processes $(\alpha_{t})_{t\leq T}$ and $(\beta_{t})_{t\leq T}$
are the intervention functions of $a_{1}$ and $a_{2}$ respectively. The agents
are also allowed to stop controlling at stopping times $\sigma$ for $a_{1}$
and $\tau$ for $a_{2}$. When one of them decides to stop controlling first,
the control of the system is stopped. The interests of the agents are
antagonistic and there is a payoff whose expression is given by:
\begin{equation}%
\begin{array}
[c]{l}%
\Gamma(\alpha,\sigma;\beta,\tau):=\mathbb{E}\left[  exp\{%
{\displaystyle\int_{0}^{T\wedge\tau\wedge\sigma}}
\varphi(s,X_{s}^{\alpha,\beta},\alpha_{s},\beta_{s})ds+h(\sigma,X_{\sigma
}^{\alpha,\beta})1_{[\sigma\leq\tau<T]}\right.  \\
\qquad\qquad\qquad\qquad\qquad\qquad\qquad\qquad\left.  +h^{\prime}%
(\tau,X_{\tau}^{\alpha,\beta})1_{[\tau<\sigma]}+g(X_{T}^{\alpha,\beta
})1_{[\sigma=\tau=T]}\}\right]  .
\end{array}
\label{payoff}%
\end{equation}
This payoff $\Gamma(\alpha,\sigma;\beta,\tau)$ is a reward for $a_{1}$ and a
cost for $a_{2}$, therefore the first agent aims at maximizing it while the
second one, his objective is to minimize the same quantity. The role of the
exponential utility function is to capture the sensitiveness with respect to
risk of the agents. To day there have been several papers which deal with
risk-sensitive control/games (see e.g. \cite{snic, bensnagai,bensnagai2, nh,
dupuis}). The risk neutral game corresponds to the case when the payoff equals
to the expectation of the quantity inside exponential in (\ref{payoff}). So
one of the objectives of this paper is to study the upper and lower values of
this mixed zero-sum stochastic differential game which are defined by:
\begin{equation}
\sup_{\alpha}\inf_{\beta}\sup_{\sigma}\inf_{\tau}\Gamma(\alpha,\sigma
;\beta,\tau)\quad\mbox{ and }\quad\inf_{\beta}\sup_{\alpha}\inf_{\tau}%
\sup_{\sigma}\Gamma(\alpha,\sigma;\beta,\tau).\label{Gamma definition}%
\end{equation}
When the data of this problem do not depend on the controls $\alpha$ and
$\beta$ and the criterion is of risk-neutral type, the problem reduces to the
well known zero-sum Dynkin game which is involved when dealing with American
game options. On the other hand, in order to obtain the Ramsey model it is
enough to take $\tau=\sigma=T$ and assume that the data do not depend on
$\beta$.

Generally speaking, in this paper we are going to consider a more general
setting of payoffs, namely payoffs defined by solutions of BSDEs with two
reflecting barriers and continuous coefficients whose growth are quadratic
with respect to the component $z$ for which the payoff (\ref{payoff}) is
somehow a particular case.

The BSDEs with two reflecting barriers have been introduced by
Cvitanic-Karatzas in \cite{CK}. They have generalized a previous work by
El-Karoui et al. \cite{elkal} on BSDEs with only one reflecting barrier. Since
then BSDEs with two reflecting barriers have attracted a lot of research
activity especially in connection with zero-sum Dynkin games and American game
options (see \cite{BHM,HH1,HH2,H}). In \cite{HH2}, the authors have shown that
if the barriers are completely separated and the coefficient of the BSDE is
continuous with quadratic growth, then a minimal and a maximal solution exist
for the BSDE. This is the framework which we adopt along with this paper.

For decades there have been a lot of research activity on stochastic
differential games (see e.g. \cite{bcr,BL1,BL2, FS, snic,hl1,hl2,hl3, hlp,
hlw},... and the references therein ). The lower and upper values of a
zero-sum stochastic differential game have been already investigated by
Fleming-Souganidis in \cite{FS}. They proved that they are unique solutions in
viscosity sense of their associated Hamilton-Jacobi-Bellman-Isaacs equations.
The upper and lower values of a zero-sum mixed differential game have been
studied in \cite{hl3,HH2,hlw} in connection with reflected BSDEs. For this
latter type of games, another framework is the one considered in \cite{BL2}.
The authors studied the mixed zero-sum differential game when the dynamics of
the controlled system is solution of (\ref{SDEintro}) and the payoffs are
given by a solution of a controlled BSDE with two reflecting barriers whose
coefficient is Lipschitz in $(y,z)$. They have shown that the values of the
game are unique viscosity solutions of their related Hamilton-Jacobi-Bellman
equations. In this work we mainly focus on the lower and upper values of the
zero-sum stochastic differential game when the payoffs are given by a solution
of a two reflecting barrier BSDEs with continuous coefficients and quadratic
growth condition $w.r.t.$ $z$. A particular case of those BSDEs is connected
with the risk-sensitive payoff defined in (\ref{payoff}). The dynamics of the
controlled system is given by (\ref{SDEintro}). We show that those values are
unique solutions of their associated Hamilton-Jacobi-Bellman equations.

This paper is organized as follows:

In Section 2, we recall the main results related to BSDEs with two reflecting
barriers, while in Section 3 and Section 4, we give the main result of the
paper. We show that the lower and upper values are unique viscosity solutions
of their respective HJB equations. We begin to show the much more involved
issue of uniqueness and later the one of existence. In Section 5, we make the
connection between the payoffs $\Gamma$ given in (\ref{payoff}) with the
payoffs defined through solutions of a specific BSDEs with two reflecting
barriers considered along with this article. $\Box$

\section{Preliminaries. Notations. Hypothesis}

The purpose of this section is to introduce some basic notations and results
concerning RBSDEs with two barriers, which will be needed throughout this
paper. In all that follows we shall consider a finite horizon $T>0$ and a
complete probability\ space $\left(  \Omega,\mathcal{F},\mathbb{P}\right)
\mathcal{\ }$on which is defined a standard $d$-dimensional Brownian motion
$W=\left(  W_{t}\right)  _{t\leq T}$ whose natural filtration is denoted
$\mathbb{F=}\{\mathcal{F}_{t},\ 0\leq t\leq T\}.$ More precisely, $\mathbb{F}$
is the filtration generated by the process $W$ and augmented by $\mathcal{N}%
_{\mathbb{P}}$, the set of all $\mathbb{P}$-null sets, \emph{i.e.}
$\mathcal{F}_{t}=\sigma\{W_{s},$\ $s\leq t\}\vee\mathcal{N}_{\mathbb{P}}%
$.$\smallskip$

Let us consider:

\begin{enumerate}
\item[(i)] $\mathcal{P}$, the $\sigma$-algebra of \thinspace$\mathcal{F}_{t}%
$-progressively measurable sets on $\left[  0,T\right]  \times\Omega$;

\item[(ii)] $\mathcal{L}^{2,k}$, the set of $\mathcal{P}$-measurable and
$\mathbb{R}^{k}$-valued processes $z=(z_{t})_{t\leq T}$ such that $\int
_{0}^{T}|z_{t}|^{2}dt<\infty$, $\mathbb{P}$-a.s.; $\mathcal{H}^{2,k}$ is the
subspace of $\mathcal{L}^{2,k}$, such that $\mathbb{E}[\int_{0}^{T}|z_{t}%
|^{2}dt]<\infty$;

\item[(iii)] $\mathcal{S}^{2}$, the set of $\mathcal{P}$-measurable and
continuous processes $Y=\left(  Y_{t}\right)  _{t\leq T}$ such that
$\mathbb{E}[\sup\limits_{t\leq T}|Y_{t}|^{2}]<\infty$;

\item[(iv)] $\mathcal{M}$, the set of continuous $\mathcal{P}$-measurable
nondecreasing processes $(K_{t})_{t\leq T}$ such that $K_{0}=0$ and
$K_{T}<\infty$, $\mathbb{P}$-a.s.\medskip
\end{enumerate}

Let us now recall the existence result for the solutions of RBSDEs with two
barriers with quadratic growth coefficient. For that, let us take four objects
which define the equation:

\begin{itemize}
\item a continuous function $F:\left[  0,T\right]  \times\Omega\times
\mathbb{R\times R}^{d}\longrightarrow\mathbb{R}$ (also known as the
\emph{coefficient or generator} of the equation) which grows subquadratically
with respect to $z$, i.e., there exists a constant $C>0$ such that,
$\mathbb{P}$-a.s.,%
\begin{equation}
\left\vert F(t,\omega,y,z)\right\vert \leq C(1+\left\vert z\right\vert
^{2}),\text{ }\forall\left(  t,y,z\right)  \in\left[  0,T\right]  \times
\Omega\times\mathbb{R\times R}^{d}; \label{crpatr}%
\end{equation}

\item a terminal value $\xi$, which is a $\mathcal{F}_{T}$-measurable random variable;

\item two processes $U=\left(  U_{t}\right)  _{t\leq T}$ and $L=\left(
L_{t}\right)  _{t\leq T}$ from $\mathcal{S}^{2}$, satisfying $L_{t}<U_{t}%
$,$\ $for all $t\leq T$, and $L_{T}\leq\xi\leq U_{T}.$
\end{itemize}

\begin{definition}
A solution for the RBSDE associated with $\left(  F,\xi,L,U\right)  $ is a
quadruple of $\mathcal{P}$-measurable processes $\left(  Y_{t},Z_{t},K_{t}%
^{+},K_{t}^{-}\right)  _{t\leq T}$ from $\mathcal{S}^{2}\times\mathcal{L}%
^{2,d}\times\mathcal{M}\times\mathcal{M}$ such that, $\mathbb{P}$-a.s.%
\begin{equation}
Y_{t}=\xi+\int_{t}^{T}F\left(  s,Y_{s},Z_{s}\right)  ds+\left(  K_{T}%
^{+}-K_{t}^{+}\right)  -\left(  K_{T}^{-}-K_{t}^{-}\right)  -\int_{t}^{T}%
Z_{s}dW_{s},\ \forall t\leq T\nonumber
\end{equation}
and%
\[
L_{t}\leq Y_{t}\leq U_{t},\ \forall t\leq T,\quad\int_{0}^{T}\left(
Y_{s}-L_{s}\right)  dK_{s}^{+}=%
{\displaystyle\int_{0}^{T}}
\left(  U_{s}-Y_{s}\right)  dK_{s}^{-}=0.
\]

\end{definition}

We will also assume that $U$, $L$ and $\xi$ are bounded, i.e.,
\begin{equation}
\operatorname*{esssup}[\left\vert \xi\right\vert +\sup_{t\leq T}\left\{
\left\vert U_{t}\right\vert +\left\vert L_{t}\right\vert \right\}  ] <+\infty.
\label{bound}%
\end{equation}
We have the following results (see, \emph{i.e}, \cite{HH2}, Theorem 3.2 and
Remark 3.3).

\begin{theorem}
Under the assumptions (\ref{crpatr}) and (\ref{bound}), there exists a
$\mathcal{P}$-measurable process $\left(  Y,Z,K^{+},K^{-}\right)  $ solution
for the RBSDE associated with $\left(  F,\xi,L,U\right)  $. Moreover, the
solution is maximal, i.e., if $\left(  Y^{\prime},Z^{\prime},K^{\prime
+},K^{\prime-}\right)  $ is another solution of the above equation, then
$\mathbb{P}$-a.s., for all $t\leq T$, we have $Y_{t}^{\prime}\leq Y_{t}$.
\end{theorem}

The main idea for proving the existence of a solution for RBSDE associated
with $\left(  F,\xi,L,U\right)  $ is to find a solution for a RBSDE associated
with data obtained by an exponential transform of $\left(  F,\xi,L,U\right)
$. Then, by a limiting procedure applied to a suitably constructed monotone
sequence of bounded, continuous generators which approximate $f$, it is shown
that the initial equation has a solution. We will use a similar technique for
the proof of our existence result for the Isaacs equation associated with
stochastic games. Another useful result is the following comparison principle.

\begin{proposition}
[comparison]Let $F$ and $F^{\prime}$ be two generators satisfying
(\ref{crpatr}) and (\ref{bound}) such that, $\mathbb{P}$-a.s., $F\left(
t,\omega,y,z\right)  \leq F^{\prime}\left(  t,\omega,y,z\right)  $ for any
$t$, $y$, $z$, and consider $\left(  Y_{t},Z_{t},K_{t}^{+},K_{t}^{-}\right)
_{t\leq T}$,\newline(resp. $\left(  Y_{t}^{\prime},Z_{t}^{\prime}%
,K_{t}^{\prime+},K_{t}^{\prime-}\right)  _{t\leq T}$) the maximal solution of
the RBSDE associated with $\left(  F,\xi,L,U\right)  $ (resp. $\left(
F^{\prime},\xi,L,U\right)  $). Then, $\mathbb{P}$-a.s. $Y_{t}\leq
Y_{t}^{\prime}$, for all $t\leq T$.
\end{proposition}

\section{Main results}

We introduce now the framework for the study of stochastic differential games
with reflection for two players and we give the results which assert the
existence and uniqueness of viscosity solutions for the associated Isaacs
equations with obstacles. The proofs of these results will be detailed in the
next section.

\subsection{The setting of the problem}

Let $A$ and $B$ be two compact metric spaces.

\begin{definition}
An admissible control process $\alpha=\left(  \alpha_{s}\right)  _{s\in\left[
t,T\right]  }$ (resp., $\beta=\left(  \beta_{s}\right)  _{s\in\left[
t,T\right]  }$) for Player I (resp., Player II) on $\left[  t,T\right]  $
($t<T$) is an $\mathbb{F}$-progressively measurable process taking values in
$A$ (resp., $B$). The set of all admissible controls on $\left[  t,T\right]  $
for the two players will be denoted by $\mathcal{A}_{t}$, respectively
$\mathcal{B}_{t}$.
\end{definition}

Now, for $t<T$, $\alpha\left(  \cdot\right)  \in\mathcal{A}_{t}$ and
$\beta\left(  \cdot\right)  \in\mathcal{B}_{t}$, let us consider the following
SDE:
\begin{equation}
\left\{
\begin{array}
[c]{l}%
dX_{s}^{t,x;\alpha,\beta}=b(s,X_{s}^{t,x;\alpha,\beta},\alpha_{s},\beta
_{s})ds+\sigma(s,X_{s}^{t,x;\alpha,\beta},\alpha_{s},\beta_{s})dW_{s}%
,\ s\in\left[  t,T\right]  ;\smallskip\\
X_{s}^{t,x;\alpha,\beta}=x,\ s\leq t,
\end{array}
\right.  \label{SDE}%
\end{equation}
where the coefficients
\[
b:\left[  0,T\right]  \times\mathbb{R}^{n}\times A\times B\longrightarrow
\mathbb{R}^{n}\text{ and }\sigma:\left[  0,T\right]  \times\mathbb{R}%
^{n}\times A\times B\longrightarrow\mathbb{R}^{n\times d}%
\]
satisfy the following condition:%
\begin{equation}
\left\{
\begin{array}
[c]{l}%
\ \left(  i\right)  \quad\text{for every }x\in\mathbb{R}^{n},\text{ }b\left(
\cdot,x,\cdot,\cdot\right)  \text{ and}\ \sigma\left(  \cdot,x,\cdot
,\cdot\right)  \text{ are continuous; moreover,}\smallskip\\
\left(  ii\right)  \quad\text{there exists }C_{L}>0\text{ such that, for all
}t\in\left[  0,T\right]  ,\ x,x^{\prime}\in\mathbb{R}^{n},\ \alpha\in
A,\ \beta\in B,\smallskip\\
\quad\quad|b\left(  t,x,\alpha,\beta\right)  -b\left(  t,x^{\prime}%
,\alpha,\beta\right)  |+|\sigma\left(  t,x,\alpha,\beta\right)  -\sigma\left(
t,x^{\prime},\alpha,\beta\right)  |\leq C_{L}|x-x^{^{\prime}}|.
\end{array}
\right.  \tag{H1}\label{H1}%
\end{equation}
It is clear (for example, see \cite{KS}) that, under the assumptions
(\ref{H1}), for every $\left(  \alpha\left(  \cdot\right)  ,\beta\left(
\cdot\right)  \right)  \in\mathcal{A}_{t}\times\mathcal{B}_{t}$, the SDE
(\ref{SDE}) has a unique solution. Moreover, for every $p\geq2$, there exists
$C_{p}>0$ such that, for all $t\in\left[  0,T\right]  $,$\ x,x^{\prime}%
\in\mathbb{R}^{n}$, $\left(  \alpha\left(  \cdot\right)  ,\beta\left(
\cdot\right)  \right)  \in\mathcal{A}_{t}\times\mathcal{B}_{t}$, we have,
a.s.:%
\[%
\begin{array}
[c]{l}%
\mathbb{E}\left[  \left.  \sup\limits_{s\in\lbrack t,T]}|X_{s}^{t,x;\alpha
,\beta}-X_{s}^{t,x^{\prime};\alpha,\beta}|^{p}\right\vert \mathcal{F}%
_{t}\right]  \leq C_{p}|x-x^{\prime}|^{p};\medskip\\
\mathbb{E}\left[  \left.  \sup\limits_{s\in\lbrack t,T]}|X_{s}^{t,x;\alpha
,\beta}|^{p}\right\vert \mathcal{F}_{t}\right]  \leq C_{p}(1+|x|^{p});
\end{array}
\]
one can see e.g. \cite{BL1, BL2} for more details. The constant $C_{p}$
depends only on the Lipschitz and the linear growth constants of $b$ and
$\sigma$.\medskip

Let us now consider the functions
\[
g:\mathbb{R}^{n}\longrightarrow\mathbb{R},\ h,h^{\prime}:\left[  0,T\right]
\times\mathbb{R}^{n}\longrightarrow\mathbb{R},\ F:\left[  0,T\right]
\times\mathbb{R}^{n}\times\mathbb{R}\times\mathbb{R}^{d}\times A\times
B\longrightarrow\mathbb{R}%
\]
that satisfy the following conditions:%
\begin{equation}
\left\{
\begin{array}
[c]{l}%
\ \left(  i\right)  \quad g\text{ is continuous and bounded; }h\text{ and
}h^{^{\prime}}\text{ are also continuous and bounded}\smallskip\\
\text{and, for any }(t,x)\in\left[  0,T\right]  \times\mathbb{R}^{n},\text{
}h(t,x)<h^{\prime}(t,x)\text{. Moreover, we assume that}\smallskip\\
\quad\quad\quad\quad\quad\quad\quad\quad h\left(  T,x\right)  \leq g\left(
x\right)  \leq h^{\prime}\left(  T,x\right)  ,\ \forall x\in\mathbb{R}%
^{n};\smallskip\\
\left(  ii\right)  \quad F\text{ is continuous and has quadratic growth in
}z,\text{ }\emph{i.e.}\smallskip\\
\left\vert F\left(  t,x,y,z,\alpha,\beta\right)  \right\vert \leq
C(1+\left\vert z\right\vert ^{2}),\text{ }\forall\left(  t,x,y,z,\alpha
,\beta\right)  \in\left[  0,T\right]  \times\mathbb{R}^{n}\times
\mathbb{R}\times\mathbb{R}^{d}\times A\times B.
\end{array}
\right.  \tag{H2}\label{bariere}%
\end{equation}
Under the above hypothesis, for any $\left(  t,x\right)  \in\left[
0,T\right]  \times\mathbb{R}^{n}$ and $\left(  \alpha\left(  \cdot\right)
,\beta\left(  \cdot\right)  \right)  \in\mathcal{A}_{t}\times\mathcal{B}_{t}$,
there exists a maximal solution $\left(  Y^{t,x;\alpha,\beta},Z^{t,x;\alpha
,\beta},K^{+,t,x;\alpha,\beta},K^{-,t,x;\alpha,\beta}\right)  $ of the RBSDE
associated with
\begin{equation}
\left(  F\left(  \cdot,X^{t,x;\alpha,\beta},\cdot,\cdot,\alpha\left(
\cdot\right)  ,\beta\left(  \cdot\right)  \right)  ,g(X_{T}^{t,x;\alpha,\beta
}),h\left(  \cdot,X^{t,x;\alpha,\beta}\right)  ,h^{\prime}\left(
\cdot,X^{t,x;\alpha,\beta}\right)  \right)  ,\label{mainRBSDE}%
\end{equation}
where $X^{t,x;\alpha,\beta}$ is the solution of equation (\ref{SDE}).

\begin{definition}
A nonanticipative strategy for Player I on $\left[  t,T\right]  $ is an
application $S_{1}:\mathcal{B}_{t}\longrightarrow\mathcal{A}_{t}$ such that,
for any $\mathbb{F}$-stopping time $\tau:\Omega\longrightarrow\left[
t,T\right]  $ and any $\beta_{1}\left(  \cdot\right)  ,\beta_{2}\left(
\cdot\right)  \in\mathcal{B}_{t}$ satisfying $\beta_{1}=\beta_{2}$ on $\left[
t,\tau\right]  $, $\mathbb{P}\otimes dt\ $a.e., we have that $S_{1}\left(
\beta_{1}\right)  =S_{1}\left(  \beta_{2}\right)  $ on $\left[  t,\tau\right]
$, $\mathbb{P}\otimes dt\ $a.e. A nonanticipative strategy on $\left[
t,T\right]  $ for the second player is a function $S_{2}:\mathcal{A}%
_{t}\longrightarrow\mathcal{B}_{t}$ defined in the same manner. We will denote
the sets of nonanticipative strategies for the two players by $\mathbb{A}_{t}%
$, respectively $\mathbb{B}_{t}$.
\end{definition}

For any given control processes $\alpha\left(  \cdot\right)  \in
\mathcal{A}_{t}$ and $\beta\left(  \cdot\right)  \in\mathcal{B}_{t}$, we
consider the associated cost functional%
\begin{equation}
J\left(  t,x;\alpha,\beta\right)  :=Y_{t}^{t,x;\alpha,\beta},\quad\left(
t,x\right)  \in\left[  0,T\right]  \times\mathbb{R}^{n} \label{cost funct}%
\end{equation}
and we define the \textit{lower value function }of the stochastic differential
game with reflection
\begin{equation}
\mathcal{W}\left(  t,x\right)  :=\operatorname*{essinf}_{S_{2}\in
\mathbb{B}_{t}}\operatorname*{esssup}_{\alpha\in\mathcal{A}_{t}}J\left(
t,x;\alpha,S_{2}\left(  \alpha\right)  \right)  \label{LVF}%
\end{equation}
and the \textit{upper value function}
\begin{equation}
\mathcal{V}\left(  t,x\right)  :=\operatorname*{esssup}_{S_{1}\in
\mathbb{A}_{t}}\operatorname*{essinf}_{\beta\in\mathcal{B}_{t}}J\left(
t,x;S_{1}\left(  \beta\right)  ,\beta\right)  . \label{UVF}%
\end{equation}

\begin{remark}
The essential infimum and the essential supremum exist and should be
understood with respect to indexed families of random variables (see the
appendix of \cite{KS}, pp. 323--325).
\end{remark}

Let us now introduce the following two Isaacs equations with obstacles%
\begin{equation}
\left\{
\begin{array}
[c]{l}%
\min\left\{  u\left(  t,x\right)  -h\left(  t,x\right)  ,\max\left\{
-\dfrac{\partial u}{\partial t}\left(  t,x\right)  -H^{-}\left(
t,x,u,Du,D^{2}u\right)  ,u\left(  t,x\right)  -h^{\prime}\left(  t,x\right)
\right\}  \right\}  =0;\\
u\left(  T,x\right)  =g\left(  x\right)  ,
\end{array}
\right.  \label{Isaac1}%
\end{equation}%
\begin{equation}
\left\{
\begin{array}
[c]{l}%
\min\left\{  v\left(  t,x\right)  -h\left(  t,x\right)  ,\max\left\{
-\dfrac{\partial v}{\partial t}\left(  t,x\right)  -H^{+}\left(
t,x,v,Dv,D^{2}v\right)  ,v\left(  t,x\right)  -h^{\prime}\left(  t,x\right)
\right\}  \right\}  =0;\\
v\left(  T,x\right)  =g\left(  x\right)  ,
\end{array}
\right.  \label{Isaac2}%
\end{equation}
associated with the Hamiltonians%
\begin{multline*}
H^{-}\left(  t,x,u,q,X\right) \\
:=\sup_{\alpha\in A}\inf_{\beta\in B}\left\{  \dfrac{1}{2}Tr\left(
\sigma\sigma^{T}\left(  t,x,\alpha,\beta\right)  X\right)  +\left\langle
b\left(  t,x,\alpha,\beta\right)  ,q\right\rangle +F\left(  t,x,u,q\sigma
\left(  t,x,\alpha,\beta\right)  ,\alpha,\beta\right)  \right\}
\end{multline*}
and%
\begin{multline*}
H^{+}\left(  t,x,u,q,X\right) \\
:=\inf_{\beta\in B}\sup_{\alpha\in A}\left\{  \dfrac{1}{2}Tr\left(
\sigma\sigma^{T}\left(  t,x,\alpha,\beta\right)  X\right)  +\left\langle
b\left(  t,x,\alpha,\beta\right)  ,q\right\rangle +F\left(  t,x,u,q\sigma
\left(  t,x,\alpha,\beta\right)  ,\alpha,\beta\right)  \right\}  ,
\end{multline*}
for all $\left(  t,x,u,q,X\right)  \in\left[  0,T\right]  \times\mathbb{R}%
^{n}\times\mathbb{R}\times\mathbb{R}^{n}\times\mathbb{S}_{n}$ ($\mathbb{S}%
_{n}$ denotes the set of symmetric $n\times n$ matrices).\bigskip

The purpose of this article is to show that, under suitable hypothesis, the
functions $\mathcal{W}$ and $\mathcal{V}$ are the unique viscosity solutions
of equations (\ref{Isaac1}), respectively (\ref{Isaac2}). This section is
dedicated to a short survey of the main results; the detailed proofs follow in
the next section.

\begin{definition}
\begin{enumerate}
\item An upper semicontinuous function $u:[0,T]\times\mathbb{R}^{n}%
\longrightarrow\mathbb{R}$ is a viscosity subsolution of equation
(\ref{Isaac1}) if $u\left(  T,x\right)  \leq g\left(  x\right)  ,\ $for every
$x\in\mathbb{R}^{n},$ and whenever $\varphi\in C^{1,2}\left(  [0,T)\times
\mathbb{R}^{n}\right)  $ and $\left(  t,x\right)  \in\lbrack0,T)\times
\mathbb{R}^{n}$ is a maximum point for $u-\varphi$, we have
\[
\min\left\{  u\left(  t,x\right)  -h\left(  t,x\right)  ,\max\left\{
-\dfrac{\partial\varphi}{\partial t}\left(  t,x\right)  -H^{-}\left(
t,x,u,D\varphi,D^{2}\varphi\right)  ,u\left(  t,x\right)  -h^{\prime}\left(
t,x\right)  \right\}  \right\}  \leq0;
\]

\item A lower semicontinuous function $u:[0,T)\times\mathbb{R}^{n}%
\longrightarrow\mathbb{R}$ is a viscosity supersolution of equation
(\ref{Isaac1}) if $u\left(  T,x\right)  \geq g\left(  x\right)  ,\ $for every
$x\in\mathbb{R}^{n},$ and whenever $\varphi\in C^{1,2}\left(  [0,T)\times
\mathbb{R}^{n}\right)  $ and $\left(  t,x\right)  \in\lbrack0,T)\times
\mathbb{R}^{n}$ is a minimum point for $u-\varphi$, we have
\[
\min\left\{  u\left(  t,x\right)  -h\left(  t,x\right)  ,\max\left\{
-\dfrac{\partial\varphi}{\partial t}\left(  t,x\right)  -H^{-}\left(
t,x,u,D\varphi,D^{2}\varphi\right)  ,u\left(  t,x\right)  -h^{\prime}\left(
t,x\right)  \right\}  \right\}  \geq0;
\]

\item A function $u:[0,T]\times\mathbb{R}^{n}\longrightarrow\mathbb{R}$ is
called a viscosity solution of equation (\ref{Isaac1}) if it is both a
viscosity sub- and supersolution for this equation.
\end{enumerate}
\end{definition}

\begin{remark}
\qquad

\begin{enumerate}
\item[(i)] Of course, the definition of viscosity solutions for equation
(\ref{Isaac2}) is similar.

\item[(ii)] In the above definitions, one can take strict local maximum or
minimum point instead of global maximum, respectively minimum point.\medskip
\end{enumerate}
\end{remark}

We recall (from \cite{CIL}, pp. 49) the definition of parabolic
\textquotedblleft superjet\textquotedblright\ and \textquotedblleft
subjet\textquotedblright\ of a function defined on a locally compact set,
notions which will be needed during the proof of the uniqueness.

\begin{definition}
For a function $u:[0,T)\times\mathbb{R}^{n}\longrightarrow\mathbb{R},$ the
second-order parabolic superjet\ of \thinspace$u$ in $\left(  t_{0}%
,x_{0}\right)  \in\lbrack0,T)\times\mathbb{R}^{n}$, denoted by $\mathcal{P}%
^{2,+}u\left(  t_{0},x_{0}\right)  $, is the set of triplets $\left(
p,q,X\right)  \in\mathbb{R\times R}^{n}\times\mathbb{S}_{n}$ satisfying, as
$\left(  t,x\right)  \rightarrow\left(  t_{0},x_{0}\right)  $%
\[
u\left(  t,x\right)  \leq u\left(  t_{0},x_{0}\right)  +p\left(
t-t_{0}\right)  +\left\langle q,x-x_{0}\right\rangle +\frac{1}{2}\left\langle
X\left(  x-x_{0}\right)  ,X\left(  x-x_{0}\right)  \right\rangle +o\left(
\left\vert t-t_{0}\right\vert +\left\vert x-x_{0}\right\vert ^{2}\right)  .
\]
Switching the inequality sign in the above relation, we get the definition of
the second-order parabolic subjet\ of \thinspace$u$ in $\left(  t_{0}%
,x_{0}\right)  $, denoted by $\mathcal{P}^{2,-}u\left(  t_{0},x_{0}\right)  $.
It is clear that $\mathcal{P}^{2,-}u=-\mathcal{P}^{2,+}\left(  -u\right)  .$
\end{definition}

One can give (see \cite{CIL}) the definition of viscosity subsolution, resp.
supersolution in terms of superjets, respectively subjets, as it follows:

\begin{proposition}
An upper semicontinuous function $u:[0,T]\times\mathbb{R}^{n}\longrightarrow
\mathbb{R}$, satisfying $u\left(  T,\cdot\right)  \leq g$, is a viscosity
subsolution of equation (\ref{Isaac1}) if and only if, for every $\left(
t,x\right)  \in\lbrack0,T)\times\mathbb{R}^{n}$ and every $\left(
p,q,X\right)  \in\mathcal{P}^{2,+}u\left(  t,x\right)  $, we have%
\[
\min\left\{  u\left(  t,x\right)  -h\left(  t,x\right)  ,\max\left\{
-p-H^{-}\left(  t,x,u(t,x),q,X\right)  ,u\left(  t,x\right)  -h^{\prime
}\left(  t,x\right)  \right\}  \right\}  \leq0.
\]

\end{proposition}

A similar result holds for viscosity supersolutions.

Since $H^{-}$ is continuous, one can replace the superjets and subjets with
their closure, whose definition is given below:

\begin{definition}
For $u:[0,T)\times\mathbb{R}^{n}\longrightarrow\mathbb{R},$ $\left(
t_{0},x_{0}\right)  \in\lbrack0,T)\times\mathbb{R}^{n},$ we define
$\mathcal{\bar{P}}^{2,+}u\left(  t_{0},x_{0}\right)  $ as the set of triplets
$\left(  p,q,X\right)  \in\mathbb{R\times R}^{n}\times\mathbb{S}_{n}$ for
which there exists a sequence $\left(  t_{n},x_{n},p_{n},q_{n},X_{n}\right)
\in\lbrack0,T)\times\mathbb{R}^{n}\times\mathbb{R\times R}^{n}\times
\mathbb{S}_{n}$ such that $\left(  p_{n},q_{n},X_{n}\right)  \in
\mathcal{P}^{2,+}u\left(  t_{n},x_{n}\right)  ,$ for all $n\in\mathbb{N}$, and%
\[
\left(  t_{n},x_{n},u\left(  t_{n},x_{n}\right)  ,p_{n},q_{n},X_{n}\right)
\rightarrow\left(  t_{0},x_{0},u\left(  t_{0},x_{0}\right)  ,p,q,X\right)  .
\]

Similarly, we define $\mathcal{\bar{P}}^{2,-}u\left(  t_{0},x_{0}\right)  $.
\end{definition}

In order to prove the uniqueness result, we will need additional properties
imposed on the generator $F$. We suppose that $F$ satisfies the following
assumption: there exist a constant $C>0$ and, for every $\varepsilon>0,$ there
exists a constant $C_{\varepsilon}>0$, such that for all $\left(
t,x,y,z,\alpha,\beta\right)  \in\left[  0,T\right]  \times\mathbb{R}^{n}%
\times\mathbb{R}\times\mathbb{R}^{d}\times A\times B$ we have%
\begin{align}
\left\vert \left(  F+\dfrac{\partial F}{\partial x}+\dfrac{\partial
F}{\partial z}\right)  \left(  t,x,y,z,\alpha,\beta\right)  \right\vert  &
\leq C(1+\left\vert z\right\vert ^{2}),\tag{H3}\label{marginiri}\\
\dfrac{\partial F}{\partial y}\left(  t,x,y,z,\alpha,\beta\right)   &  \leq
C_{\varepsilon}+\varepsilon\left\vert z\right\vert ^{2}.\nonumber
\end{align}

\subsection{Results}

The framework is now set for the main part of this paper.

\begin{theorem}
[Existence]Under the assumptions (\ref{H1}), (\ref{bariere}), and
(\ref{marginiri}), the \textit{lower value function }$\mathcal{W}$ defined by
(\ref{LVF}) is a viscosity solution of the Isaacs equation with two barriers
(\ref{Isaac1}), while the \textit{upper value function }$\mathcal{V}$ defined
by (\ref{UVF}) is a viscosity solution of the Isaacs equation (\ref{Isaac2}).
\end{theorem}

\begin{theorem}
[Uniqueness]\label{uniqueness}Under the assumptions (\ref{H1}), (\ref{bariere}%
) and (\ref{marginiri}), if $u$ is a bounded viscosity subsolution and $v$ is
a bounded viscosity supersolution of equation (\ref{Isaac1}), then
\[
u\left(  t,x\right)  \leq v\left(  t,x\right)  ,\text{ }\forall\left(
t,x\right)  \in\left[  0,T\right]  \times\mathbb{R}^{n}.
\]

The same comparison principle holds for the Isaacs equation (\ref{Isaac2}).
\end{theorem}

\begin{remark}
If, in addition, the Isaacs' condition holds, \emph{i.e.}%
\[
H^{-}\left(  t,x,u,q,X\right)  =H^{+}\left(  t,x,u,q,X\right)  ,
\]
for every $\left(  t,x,u,q,X\right)  \in\left[  0,T\right]  \times
\mathbb{R}^{n}\times\mathbb{R}\times\mathbb{R}^{n}\times\mathbb{S}_{n}$, then
the two Isaacs equations coincide and it follows that the upper and the lower
value functions are equal, which means that the corresponding reflected
stochastic differential game has a value.
\end{remark}

\section{Proofs}

We will focus our attention on the first Isaacs equation (\ref{Isaac1}), the
case of the equation (\ref{Isaac2}) being treated in a similar manner.

\subsection{Uniqueness}

Let us consider a subsolution $u$, respectively a supersolution $v$, of
equation (\ref{Isaac1}).

First, we will make a change of variable, which preserves viscosity sub- and
supersolutions and which transforms the equation into a Isaacs equation whose
Hamiltonian will satisfy some kind of monotonicity.

We consider $\tilde{C}:=\max\left(  \left\Vert u\right\Vert _{\infty
},\left\Vert v\right\Vert _{\infty}\right)  +1$ and introduce the positive,
increasing function $\rho$ used in \cite{K}, pp. 583:%
\[
\rho:\mathbb{R\longrightarrow}\left(  -\left(  \ln\gamma\right)
/\lambda,+\infty\right)  ,\quad\rho\left(  x\right)  :=\frac{1}{\lambda}%
\ln\left(  \frac{e^{\lambda\gamma x}+1}{\gamma}\right)  ,
\]
for$\,$positive $\gamma$ and $\lambda$ satisfying $-\left(  \ln\gamma\right)
/\lambda\leq\tilde{C}$. We make the change of variable $\bar{u}:=\rho
^{-1}(e^{Kt}(u-\tilde{C}))$, with $K>0$. Equation (\ref{Isaac1}) becomes%
\begin{equation}
\left\{
\begin{array}
[c]{l}%
\min\left\{  \rho\left(  \bar{u}\right)  -\rho\left(  \bar{h}\right)
,\max\left\{  \rho^{\prime}\left(  \bar{u}\right)  \left[  -\dfrac
{\partial\bar{u}}{\partial t}-\bar{H}^{-}\left(  t,x,\bar{u},D\bar{u}%
,D^{2}\bar{u}\right)  \right]  ,\rho\left(  \bar{u}\right)  -\rho\left(
\bar{h}^{\prime}\right)  \right\}  \right\}  =0;\\
\bar{u}\left(  T,x\right)  =\bar{g}\left(  x\right)  ,
\end{array}
\right.  \label{pb transformata}%
\end{equation}
where%
\[
\bar{H}^{-}\left(  t,x,\bar{u},\bar{q},\bar{X}\right)  =\sup_{\alpha\in A}%
\inf_{\beta\in B}\left[  \frac{1}{2}\operatorname*{Tr}\left(  \sigma\sigma
^{T}\left(  t,x,\alpha,\beta\right)  \bar{X}\right)  +\left\langle b\left(
t,x,\alpha,\beta\right)  ,\bar{q}\right\rangle +\bar{F}\left(  t,x,\bar
{u},\bar{q}\sigma\left(  t,x,\alpha,\beta\right)  ,\alpha,\beta\right)
\right]  ,
\]
with $\bar{F}$ defined by%
\[
\bar{F}\left(  t,x,\bar{u},\bar{z},\alpha,\beta\right)  =\dfrac{\rho
^{\prime\prime}\left(  \bar{u}\right)  }{\rho^{\prime}\left(  \bar{u}\right)
}\left\vert \bar{z}\right\vert ^{2}-K\dfrac{\rho\left(  \bar{u}\right)  }%
{\rho^{\prime}\left(  \bar{u}\right)  }+\frac{e^{Kt}}{\rho^{\prime}\left(
\bar{u}\right)  }F(t,x,e^{-Kt}\rho\left(  \bar{u}\right)  +\tilde{C}%
,e^{-Kt}\rho^{\prime}\left(  \bar{u}\right)  \bar{z},\alpha,\beta)
\]
and%
\begin{align*}
\bar{h}\left(  t,x\right)   &  :=\rho^{-1}(e^{Kt}(h\left(  t,x\right)
-\tilde{C})),\\
\bar{h}^{\prime}\left(  t,x\right)   &  :=\rho^{-1}(e^{Kt}(h^{\prime}\left(
t,x\right)  -\tilde{C})),\\
\bar{g}\left(  x\right)   &  :=\rho^{-1}(e^{Kt}(g\left(  x\right)  -\tilde
{C})).
\end{align*}
The function $\bar{F}$ verifies, for $\gamma$ big enough,

\begin{lemma}
There exist some positive constants $\tilde{K}$ and $\bar{C}$ such that for
all $t\in\left(  0,T\right)  ,\ \left(  \alpha,\beta\right)  \in A\times B$,
$x,y\in\mathbb{R}^{n},\ z,z^{\prime}\in\mathbb{R}^{d}$, and $u,v\in\mathbb{R}$
such that $u<v$,%
\[
\bar{F}\left(  t,x,u,z,\alpha,\beta\right)  -\bar{F}\left(  t,y,v,z^{\prime
},\alpha,\beta\right)  \leq\mathcal{K}\left(  z,z^{\prime}\right)  (-\tilde
{K}\left(  u-v\right)  +\bar{C}|x-y|+\bar{C}|z-z^{\prime}|),
\]
where $\mathcal{K}\left(  z,z^{\prime}\right)  :=\left(  1+\dfrac{|z|^{2}}%
{2}+\dfrac{|z^{\prime}|^{2}}{2}\right)  .$
\end{lemma}

As we said before, by this transformation, if $u$ (resp., $v$) is a
subsolution (resp., a supersolution) of equation (\ref{Isaac1}), then $\bar
{u}$ (resp., $\bar{v}$) is one for equation (\ref{pb transformata}).
Therefore, we want to prove that%
\[
M:=\sup_{\left(  t,x\right)  \in\left[  0,T\right]  \times\mathbb{R}^{n}}%
(\bar{u}\left(  t,x\right)  -\bar{v}\left(  t,x\right)  )
\]
is negative. Define also
\[
M\left(  h\right)  :=\sup_{\left\vert x-y\right\vert \leq h}|\bar{u}\left(
t,x\right)  -\bar{v}\left(  t,y\right)  |\quad\text{and}\quad M^{\prime}%
:=\lim_{h\rightarrow0}M\left(  h\right)  .
\]
It is clear that $M\leq M^{\prime}$.$\smallskip$

We assume to the contrary $M>0$ and define for every $\varepsilon,$ $\eta>0$%
\[
\Psi_{\varepsilon,\eta}\left(  t,x,y\right)  :=\bar{u}\left(  t,x\right)
-\bar{v}\left(  t,y\right)  -\dfrac{|x-y|^{2}}{\varepsilon^{2}}-\eta
(\left\vert x\right\vert ^{2}+\left\vert y\right\vert ^{2}).
\]
Let us consider
\[
M_{\varepsilon,\eta}:=\sup\limits_{\left(  t,x,y\right)  \in\left[
0,T\right]  \times\mathbb{R}^{n}\times\mathbb{R}^{n}}\Psi_{\varepsilon,\eta
}\left(  t,x,y\right)  =\max\limits_{\left(  t,x,y\right)  \in\left[
0,T\right]  \times\mathbb{R}^{n}\times\mathbb{R}^{n}}\Psi_{\varepsilon,\eta
}\left(  t,x,y\right)
\]
Since the functions $\bar{u}$ and $\bar{v}$ are bounded, the supremum of
$\Psi_{\varepsilon,\eta}$ is reached at some point $\left(  t^{\varepsilon
,\eta},x^{\varepsilon,\eta},y^{\varepsilon,\eta}\right)  $, which will be
denoted for simplicity $\left(  \hat{t},\hat{x},\hat{y}\right)  .$ We will use
this notations each time we do not want to make explicit the dependence on
$\varepsilon$ and $\eta$.

Let recall some notations and results from \cite{K}, pp. 586-587: for the
sequence $\left(  a_{\varepsilon,\eta}\right)  _{\varepsilon,\eta}$, we write
$a=\lim\limits_{\varepsilon\ll\eta\rightarrow0}a_{\varepsilon,\eta}$ if%
\[
\lim\inf\nolimits_{\varepsilon\ll\eta\rightarrow0}a_{\varepsilon,\eta}%
=\lim\sup\nolimits_{\varepsilon\ll\eta\rightarrow0}a_{\varepsilon,\eta}=a,
\]
where%
\[%
\begin{array}
[c]{l}%
\lim\inf\nolimits_{\varepsilon\ll\eta\rightarrow0}a_{\varepsilon,\eta}%
=\lim\inf\nolimits_{\eta\rightarrow0}\left(  \lim\inf\nolimits_{\varepsilon
\rightarrow0}a_{\varepsilon,\eta}\right)  ,\medskip\\
\lim\sup\nolimits_{\varepsilon\ll\eta\rightarrow0}a_{\varepsilon,\eta}%
=\lim\sup\nolimits_{\eta\rightarrow0}\left(  \lim\sup\nolimits_{\varepsilon
\rightarrow0}a_{\varepsilon,\eta}\right)  .\medskip
\end{array}
\]
The following result is the equivalent of Lemma 3.1 of \cite{CIL}.

\begin{lemma}
Considering the above notations, we have%
\[%
\begin{array}
[c]{ll}%
\left(  i\right)   & \lim\limits_{\varepsilon\ll\eta\rightarrow0}%
M_{\varepsilon,\eta}=M,\quad\lim\limits_{\varepsilon\ll\eta\rightarrow0}%
\bar{u}\left(  \hat{t},\hat{x}\right)  -\bar{v}\left(  \hat{t},\hat{y}\right)
=M;\smallskip\\
\left(  ii\right)   & \lim\limits_{\eta\ll\varepsilon\rightarrow
0}M_{\varepsilon,\eta}=M^{\prime},\quad\lim\limits_{\eta\ll\varepsilon
\rightarrow0}\bar{u}\left(  \hat{t},\hat{x}\right)  -\bar{v}\left(  \hat
{t},\hat{y}\right)  =M^{\prime};\smallskip\\
\left(  iii\right)   & \lim\limits_{\eta\ll\varepsilon\rightarrow0}%
\dfrac{\left\vert \hat{x}-\hat{y}\right\vert }{\varepsilon}=0,\quad
\lim\limits_{\eta\ll\varepsilon\rightarrow0}\eta(\left\vert \hat{x}\right\vert
^{2}+\left\vert \hat{y}\right\vert ^{2})=0.\smallskip
\end{array}
\]

\end{lemma}

So, by extracting a subsequence, we suppose that for every $\eta,$ the
sequence $\left(  t^{\varepsilon,\eta}\right)  _{\varepsilon}$ converges to a
limit $t^{\eta}$ as $\varepsilon\rightarrow0$ and, extracting again a
subsequence, the sequences $\left(  x^{\varepsilon,\eta}\right)
_{\varepsilon}$ and $\left(  y^{\varepsilon,\eta}\right)  _{\varepsilon}$
converge to a common limit $x^{\eta}.$

Now define the functions $\phi_{1},$ $\phi_{2}:\left[  0,T\right]
\times\mathbb{R}^{n}\longrightarrow\mathbb{R}$ by%
\begin{align*}
\phi_{1}\left(  t,x\right)   &  =\bar{v}\left(  t,y^{\varepsilon,\eta}\right)
+\dfrac{|x-y^{\varepsilon,\eta}|^{2}}{\varepsilon^{2}}+\eta(\left\vert
x\right\vert ^{2}+\left\vert y^{\varepsilon,\eta}\right\vert ^{2})\\
\phi_{2}\left(  t,y\right)   &  =\bar{u}\left(  t,x^{\varepsilon,\eta}\right)
-\dfrac{|x^{\varepsilon,\eta}-y|^{2}}{\varepsilon^{2}}-\eta(\left\vert
x^{\varepsilon,\eta}\right\vert ^{2}+\left\vert y\right\vert ^{2}).
\end{align*}
It is obvious that $\left(  t^{\varepsilon,\eta},x^{\varepsilon,\eta}\right)
$ is a maximum point for the function $\left(  t,x\right)  \longmapsto
\Psi_{\varepsilon,\eta}\left(  t,x,y^{\varepsilon,\eta}\right)  =\left(
\bar{u}-\phi_{1}\right)  \left(  t,x\right)  $, while $\left(  t^{\varepsilon
,\eta},y^{\varepsilon,\eta}\right)  $ is a minimum point for $\left(
t,y\right)  \longmapsto-\Psi_{\varepsilon,\eta}\left(  t,x^{\varepsilon,\eta
},y\right)  =\left(  \phi_{2}-\bar{v}\right)  \left(  t,y\right)  .$ Since
$\bar{u}$ and $\bar{v}$ are viscosity subsolution, respectively supersolution,
we obtain

\begin{itemize}
\item either $t^{\varepsilon,\eta}=T$ and then $\bar{u}\left(
T,x^{\varepsilon,\eta}\right)  \leq\bar{g}\left(  x^{\varepsilon,\eta}\right)
$ and $\bar{g}\left(  y^{\varepsilon,\eta}\right)  \leq\bar{v}\left(
T,y^{\varepsilon,\eta}\right)  $,

\item or $t^{\varepsilon,\eta}\neq T$ and we have, in $\left(  \hat{t},\hat
{x}\right)  $,%
\[
\min\left\{  \rho(\bar{u})-\rho(\bar{h}),\max\left\{  \rho^{\prime}(\bar
{u})\left[  -\dfrac{\partial\phi_{1}}{\partial t}-\bar{H}^{-}\left(  \hat
{t},\hat{x},\bar{u},D\phi_{1},D^{2}\phi_{1}\right)  \right]  ,\rho(\bar
{u})-\rho(\bar{h}^{\prime})\right\}  \right\}  \leq0
\]
and, respectively, in $\left(  \hat{t},\hat{y}\right)  $,%
\[
\min\left\{  \rho(\bar{v})-\rho(\bar{h}),\max\left\{  \rho^{\prime}(\bar
{v})\left[  -\dfrac{\partial\phi_{2}}{\partial t}-\bar{H}^{-}\left(  \hat
{t},\hat{y},\bar{v},D\phi_{2},D^{2}\phi_{2}\right)  \right]  ,\rho(\bar
{v})-\rho(\bar{h}^{\prime})\right\}  \right\}  \geq0.
\]

\end{itemize}

In the first situation, there exists a subsequence of $\left(  t^{\eta
}\right)  _{\eta}$, supposed, without restricting the generality, to be the
same, such that $t^{\eta}=T$ for every $\eta$. The semicontinuity of the
functions $\bar{u}$ and $\bar{v}$ and the continuity of $\bar{g}$ give us
that, for all $\eta$ and $\varepsilon$ sufficiently small,%
\[
\bar{u}(\hat{t},\hat{x})\leq\bar{u}(T,x^{\eta})+\eta\leq\bar{g}(x^{\eta}%
)+\eta\quad\text{ and }\quad\bar{g}(x^{\eta})-\eta\leq\bar{v}(T,x^{\eta}%
)-\eta\leq\bar{v}(\hat{t},\hat{x}).
\]
So, $\bar{u}\left(  \hat{t},\hat{x}\right)  \leq\bar{v}\left(  \hat{t},\hat
{x}\right)  +2\eta$ and, from here, we find that $M_{\varepsilon,\eta}%
\leq2\eta.$ Taking the limit of $\varepsilon$, and after this $\eta$, to zero,
it follows that $M\leq0$, which is a contradiction. Therefore $t^{\varepsilon
,\eta}\neq T$.

Now, since $\rho$ is an increasing function, if there exists a subsequence
$t^{\eta}\neq T$, and then, for every $\eta$ a subsequence of $\left(
x^{\varepsilon,\eta}\right)  _{\varepsilon}$ and one of $\left(
y^{\varepsilon,\eta}\right)  _{\varepsilon}$ such that the following
inequalities hold%
\[
\bar{u}(t^{\varepsilon,\eta},x^{\varepsilon,\eta})-\bar{h}(t^{\varepsilon
,\eta},x^{\varepsilon,\eta})\leq0\quad\text{and}\quad\bar{v}(t^{\varepsilon
,\eta},y^{\varepsilon,\eta})\geq\bar{h}(t^{\varepsilon,\eta},y^{\varepsilon
,\eta}),
\]
we obtain%
\[
M_{\varepsilon,\eta}\leq\bar{u}(t^{\varepsilon,\eta},x^{\varepsilon,\eta
})-\bar{v}(t^{\varepsilon,\eta},y^{\varepsilon,\eta})\leq\bar{h}%
(t^{\varepsilon,\eta},x^{\varepsilon,\eta})-\bar{h}(t^{\varepsilon,\eta
},y^{\varepsilon,\eta}).
\]
The continuity of $\bar{h}$ and the above Lemma leads again to a
contradiction. When, also on a subsequence like in the previous case,
\thinspace$\bar{u}(\hat{t},\hat{x})\leq\bar{h}^{\prime}(\hat{t},\hat{x})$,%
\[
-\dfrac{\partial\phi_{1}}{\partial t}\left(  \hat{t},\hat{x}\right)  -\bar
{H}^{-}\left(  \hat{t},\hat{x},\bar{u}\left(  \hat{t},\hat{x}\right)
,D\phi_{1}\left(  \hat{t},\hat{x}\right)  ,D^{2}\phi_{1}\left(  \hat{t}%
,\hat{x}\right)  \right)  \leq0
\]
and $\bar{v}(\hat{t},\hat{y})\geq\bar{h}^{\prime}(\hat{t},\hat{y})$, we obtain
also that $M\leq0$.

For the remaining situation, we first use Theorem 8.3 from \cite{CIL}.

Setting $\varphi\left(  x,y\right)  :=\dfrac{|x-y|^{2}}{\varepsilon^{2}}%
+\eta(\left\vert x\right\vert ^{2}+\left\vert y\right\vert ^{2})$, from the
mentioned result we have that there exist the matrices $\bar{X},\bar{Y}%
\in\mathbb{S}_{n}$ such that%
\[
\left\{
\begin{array}
[c]{l}%
(0,D_{x}\varphi\left(  \hat{x},\hat{y}\right)  ,\bar{X})\in\mathcal{\bar{P}%
}^{2,+}\bar{u}(\hat{t},\hat{x}),\medskip\\
(0,D_{y}\varphi\left(  \hat{x},\hat{y}\right)  ,-\bar{Y})\in\mathcal{\bar{P}%
}^{2,+}\left(  -\bar{v}\right)  (\hat{t},\hat{y})=-\mathcal{\bar{P}}^{2,-}%
\bar{v}(\hat{t},\hat{y}).
\end{array}
\right.  \medskip
\]
So $\left(  0,-D_{y}\varphi\left(  \hat{t},\hat{x},\hat{y}\right)  ,\bar
{Y}\right)  \in\mathcal{\bar{P}}^{2,-}\bar{v}(\hat{t},\hat{y}).\smallskip$

Moreover,
\begin{equation}
\left(
\begin{array}
[c]{cc}%
\bar{X} & 0\\
0 & -\bar{Y}%
\end{array}
\right)  \leq\dfrac{2}{\varepsilon^{2}}\left(
\begin{array}
[c]{cc}%
I & -I\\
-I & I
\end{array}
\right)  +2\eta\left(
\begin{array}
[c]{cc}%
I & 0\\
0 & I
\end{array}
\right)  .\label{ineg matr}%
\end{equation}
We have then from the definition of the sub- and superjet that%
\[
\left\{
\begin{array}
[c]{l}%
\bar{H}^{-}\left(  \hat{t},\hat{x},\bar{u}\left(  \hat{t},\hat{x}\right)
,D_{x}\varphi\left(  \hat{x},\hat{y}\right)  ,\bar{X}\right)  \geq0\text{
and}\\
\bar{H}^{-}\left(  \hat{t},\hat{y},\bar{v}\left(  \hat{t},\hat{y}\right)
,-D_{y}\varphi\left(  \hat{x},\hat{y}\right)  ,\bar{Y}\right)  \leq0,
\end{array}
\right.
\]
which implies%
\[
\bar{H}^{-}\left(  \hat{t},\hat{x},\bar{u}(\hat{t},\hat{x}),\frac{2(\hat
{x}-\hat{y})}{\varepsilon^{2}}+2\eta\hat{x},\bar{X}\right)  \geq\bar{H}%
^{-}\left(  \hat{t},\hat{y},\bar{v}(\hat{t},\hat{y}),\frac{2(\hat{x}-\hat{y}%
)}{\varepsilon^{2}}-2\eta\hat{y},\bar{Y}\right)  .
\]
By denoting
\begin{multline*}
\mathcal{\hat{K}}_{x}^{\alpha,\beta}:=\frac{1}{2}Tr(\sigma\sigma^{T}(\hat
{t},\hat{x},\alpha,\beta)\bar{X})\\
+\left\langle b\left(  \hat{t},\hat{x},\alpha,\beta\right)  ,\frac{2\left(
\hat{x}-\hat{y}\right)  }{\varepsilon^{2}}+2\eta\hat{x}\right\rangle +\bar
{F}\left(  \hat{t},\hat{x},\bar{u}\left(  \hat{t},\hat{x}\right)  ,\left(
\frac{2\left(  \hat{x}-\hat{y}\right)  }{\varepsilon^{2}}+2\eta\hat{x}\right)
\sigma\left(  \hat{t},\hat{x},\alpha,\beta\right)  ,\alpha,\beta\right)
\end{multline*}
and%
\begin{multline*}
\mathcal{\hat{K}}_{y}^{\alpha,\beta}:=\frac{1}{2}Tr(\sigma\sigma^{T}(\hat
{t},\hat{y},\alpha,\beta)\bar{Y})\\
+\left\langle b\left(  \hat{t},\hat{y},\alpha,\beta\right)  ,\frac{2\left(
\hat{x}-\hat{y}\right)  }{\varepsilon^{2}}-2\eta\hat{y}\right\rangle +\bar
{F}\left(  \hat{t},\hat{y},\bar{v}\left(  \hat{t},\hat{y}\right)  ,\left(
\frac{2\left(  \hat{x}-\hat{y}\right)  }{\varepsilon^{2}}-2\eta\hat{y}\right)
\sigma\left(  \hat{t},\hat{y},\alpha,\beta\right)  ,\alpha,\beta\right)
\end{multline*}
this can be written as
\begin{equation}
\sup_{\alpha\in A}\inf_{\beta\in B}\mathcal{\hat{K}}_{x}^{\alpha,\beta}%
\geq\sup_{\alpha\in A}\inf_{\beta\in B}\mathcal{\hat{K}}_{y}^{\alpha,\beta
}.\label{ineq_K}%
\end{equation}

Also, by denoting%
\[
\left\{
\begin{array}
[c]{ll}%
\hat{b}_{x}:=b(\hat{t},\hat{x},\alpha,\beta), & \quad\hat{b}_{y}:=b(\hat
{t},\hat{y},\alpha,\beta),\\
\hat{\sigma}_{x}:=\sigma(\hat{t},\hat{x},\alpha,\beta), & \quad\hat{\sigma
}_{y}:=\sigma(\hat{t},\hat{y},\alpha,\beta),\\
\bar{q}_{x}:=\dfrac{2(\hat{x}-\hat{y})}{\varepsilon^{2}}+2\eta\hat{x}, &
\quad\bar{q}_{y}:=\dfrac{2(\hat{x}-\hat{y})}{\varepsilon^{2}}-2\eta\hat{y},
\end{array}
\right.
\]

the inequality (\ref{ineg matr}) infers that%
\[
\frac{1}{2}Tr(\hat{\sigma}_{x}\hat{\sigma}_{x}^{T}\bar{X})\leq\frac{1}%
{2}Tr(\hat{\sigma}_{y}\hat{\sigma}_{y}^{T}\bar{Y})+\frac{1}{\varepsilon^{2}%
}|\hat{\sigma}_{x}-\hat{\sigma}_{y}|^{2}+\eta(\left\vert \hat{\sigma}%
_{x}\right\vert ^{2}+\left\vert \hat{\sigma}_{y}\right\vert ^{2})
\]
On the other hand%
\[
\left\langle \hat{b}_{x},\bar{q}_{x}\right\rangle =\left\langle \hat{b}%
_{y},\bar{q}_{y}\right\rangle +\left\langle \hat{b}_{x}-\hat{b}_{y}%
,\frac{2\left(  \hat{x}-\hat{y}\right)  }{\varepsilon^{2}}\right\rangle
+\left\langle \hat{b}_{y},2\eta\hat{y}\right\rangle +\left\langle \hat{b}%
_{x},2\eta\hat{x}\right\rangle
\]
and%
\begin{multline*}
\bar{F}\left(  \hat{t},\hat{x},\bar{u}\left(  \hat{t},\hat{x}\right)  ,\bar
{q}_{x}\hat{\sigma}_{x},\alpha,\beta\right)  \leq\bar{F}\left(  \hat{t}%
,\hat{y},\bar{v}\left(  \hat{t},\hat{y}\right)  ,\bar{q}_{y}\hat{\sigma}%
_{y},\alpha,\beta\right)  \\
+\mathcal{K}\left(  \bar{q}_{x}\hat{\sigma}_{x},\bar{q}_{y}\hat{\sigma}%
_{y}\right)  (-\tilde{K}\left(  \bar{u}\left(  \hat{t},\hat{x}\right)
-\bar{v}\left(  \hat{t},\hat{y}\right)  \right)  +\bar{C}\left\vert \hat
{x}-\hat{y}\right\vert +\bar{C}\left\vert \bar{q}_{x}\hat{\sigma}_{x}-\bar
{q}_{y}\hat{\sigma}_{y}\right\vert ).
\end{multline*}
Adding the last three relations, we obtain, for all $\left(  \alpha
,\beta\right)  \in A\times B$%
\begin{align*}
\mathcal{K}_{x}^{\alpha,\beta} &  \leq\mathcal{K}_{y}^{\alpha,\beta}%
+\frac{C_{L}^{2}}{\varepsilon^{2}}|\hat{x}-\hat{y}|^{2}+C\eta(1+\left\vert
\hat{x}\right\vert ^{2}+\left\vert \hat{y}\right\vert ^{2})+\frac{2C_{L}%
}{\varepsilon^{2}}|\hat{x}-\hat{y}|^{2}\\
&  +\mathcal{K}\left(  \bar{q}_{x}\hat{\sigma}_{x},\bar{q}_{y}\hat{\sigma}%
_{y}\right)  \left(  -\tilde{K}M_{\varepsilon,\eta}+\bar{C}|\hat{x}-\hat
{y}|+2\bar{C}\left(  C_{L}\frac{|\hat{x}-\hat{y}|^{2}}{\varepsilon^{2}}%
+\eta(\left\vert \hat{\sigma}_{x}\hat{x}\right\vert +\left\vert \hat{\sigma
}_{y}\hat{y}\right\vert )\right)  \right)  .
\end{align*}
Taking the $\sup_{\alpha\in A}\inf_{\beta\in B}$, and then passing to the
limit as $\varepsilon\rightarrow0$ and $\eta\rightarrow0$ (in this order), we
obtain%
\[
\sup_{\alpha\in A}\inf_{\beta\in B}\mathcal{\hat{K}}_{x}^{\alpha,\beta}%
\leq\sup_{\alpha\in A}\inf_{\beta\in B}\mathcal{\hat{K}}_{y}^{\alpha,\beta
}-\tilde{K}M,
\]
which contradicts (\ref{ineq_K}) if $M$ is strictly positive as we supposed.
So $M$ must be less or equal to zero and we obtain the desired comparison result.

\subsection{Existence}

As indicated in preliminaries, we make the following transform:%
\[
f\left(  t,x,y,z,\alpha,\beta\right)  :=\left\{
\begin{array}
[c]{ll}%
2Cy\left[  F\left(  t,x,\dfrac{\ln y}{2C},\dfrac{z}{2Cy},\alpha,\beta\right)
-\dfrac{\left\vert z\right\vert ^{2}}{4Cy^{2}}\right]  , & y>0\\
0, & y\leq0
\end{array}
\right.
\]
for $\left(  \alpha,\beta\right)  \in A\times B\ $and $\left(  t,\omega
,y,z\right)  \in\left[  0,T\right]  \times\Omega\times\mathbb{R}%
\times\mathbb{R}^{d}$, where $C$ is the constant from (\ref{marginiri}). We
also set, for the simplicity of notations,%
\[
\bar{h}:=\exp\left(  2Ch\right)  ,\text{ }\bar{h}^{\prime}:=\exp\left(
2Ch^{\prime}\right)  ,\text{ }\bar{g}:=\exp\left(  2Cg\right)  \text{ and
}M:=\max\left\{  \inf\left(  1/\bar{h}\right)  ,\sup\bar{h}^{\prime}\right\}
.
\]

As showed in \cite{HH2}, under conditions (\ref{H1}), (\ref{bariere}) and
(\ref{marginiri}), the RBSDE associated with
\[
(f,\bar{g}(X_{T}^{t,x;\alpha,\beta}),\bar{h}(\cdot,X^{t,x;\alpha,\beta}%
),\bar{h}^{\prime}(\cdot,X^{t,x;\alpha,\beta}))
\]
has a maximal solution, denoted $\left(  y^{t,x;\alpha,\beta},z^{t,x;\alpha
,\beta},k^{+,t,x;\alpha,\beta},k^{-,t,x;\alpha,\beta}\right)  $. In order to
use already known results on the lower value function associated with a RBSDE
with double barrier with a Lipschitz generator (see \cite{BL2}), we
approximate, in a monotone manner, the function $f$ with Lipschitz functions.
For that, we use the method considered in \cite{K}, which we describe in the sequel.

We will make use, for $p\geq0$, of functions $\rho_{p}\in C^{\infty}\left(
\mathbb{R}\right)  $ satisfying $\rho_{p}=1$ on $[-p,p]$ and $\rho_{p}=0$ on
$[-\left(  p+1\right)  ,p+1]^{\operatorname*{c}}$. Let us consider the
function%
\[
\tilde{f}\left(  t,x,y,z,\alpha,\beta\right)  :=f\left(  t,x,y,z,\alpha
,\beta\right)  \rho_{M}\left(  \dfrac{\ln y}{2C}\right)  ,
\]
which is bounded by a constant $C^{\prime}>0$. Let, for $p\in\mathbb{N}$,%
\[
\tilde{f}^{p}\left(  t,x,y,z,\alpha,\beta\right)  :=\tilde{f}\left(
t,x,y,z,\alpha,\beta\right)  \rho_{p}\left(  \left\vert x\right\vert
+\left\vert z\right\vert \right)  +\dfrac{3}{2^{p+2}};
\]
it is clear that the function $\tilde{f}^{p}$ is bounded and uniformly continuous.

Now, by a standard \textit{mollification} procedure, we approximate $\tilde
{f}^{p}$ by Lipschitz functions. Let $\theta\in C^{\infty}\left(
\mathbb{R}^{n+1+d}\right)  $ satisfying $\theta\geq0$, $\operatorname*{supp}%
\theta\subseteq B\left(  0;1\right)  $ and $\int_{\mathbb{R}^{n+1+d}}%
\theta\left(  a\right)  da=1$; we set%
\[
\tilde{f}_{\varepsilon}^{p}\left(  t,x,y,z,\alpha,\beta\right)  :=\dfrac
{1}{\varepsilon^{n+1+d}}\int_{\mathbb{R}^{n+1+d}}\theta\left(  \dfrac
{x-x^{\prime}}{\varepsilon},\dfrac{y-y^{\prime}}{\varepsilon},\dfrac
{z-z^{\prime}}{\varepsilon}\right)  \tilde{f}^{p}\left(  t,x,y,z,\alpha
,\beta\right)  dx^{\prime}dy^{\prime}dz^{\prime}.
\]
Then $\tilde{f}_{\varepsilon}^{p}$ is Lipschitz in $\left(  x,y,z\right)
\in\mathbb{R}^{n+1+d}$ and%
\[
|\tilde{f}_{\varepsilon}^{p}-\tilde{f}_{p}|\leq\eta_{\tilde{f}_{p}}\left(
\varepsilon\right)  ,
\]
where $\eta_{\tilde{f}^{p}}$ is the modulus of uniform continuity of
$\tilde{f}^{p}$. Hence, one can extract a sequence $\varepsilon_{p}\searrow0$
such that%
\[
|\tilde{f}_{\varepsilon_{p}}^{p}-\tilde{f}_{p}|\leq2^{-\left(  p+2\right)
},\ \forall p\in\mathbb{N}^{\ast}.
\]
An easy calculus shows us that $\tilde{f}_{\varepsilon_{p}}^{p}$ is upper
bounded by $C^{\prime}+2^{-p}$ and%
\[
(\tilde{f}_{\varepsilon_{p+1}}^{p+1}-\tilde{f}_{\varepsilon_{p}}^{p})\left(
t,x,y,z,\alpha,\beta\right)  \leq\tilde{f}\left(  t,x,y,z,\alpha,\beta\right)
\left(  \rho_{p+1}-\rho_{p}\right)  \left(  \left\vert x\right\vert
+\left\vert z\right\vert \right)  .
\]
If we set, for $p\in\mathbb{N}^{\ast}$,%
\[
f^{p}\left(  t,x,y,z,\alpha,\beta\right)  :=\rho_{p-1}\left(  \left\vert
x\right\vert +\left\vert z\right\vert \right)  \tilde{f}_{\varepsilon_{p}}%
^{p}\left(  t,x,y,z,\alpha,\beta\right)  +\left(  1-\rho_{p-1}\left(
\left\vert x\right\vert +\left\vert z\right\vert \right)  \right)  \left(
C^{\prime}+2^{-p}\right)  ,
\]
then the functions $f^{p}$ are still Lipschitz in $\left(  x,y,z\right)
\in\mathbb{R}^{n+1+d}$ and $f^{p}\searrow\tilde{f}$, the convergence being
uniform. \medskip

Consider now, for $p\in\mathbb{N}^{\ast}$, the RBSDE associated with
$(f^{p},\bar{g}(X_{T}^{t,x;\alpha,\beta}),\bar{h}(\cdot,X^{t,x;\alpha,\beta
}),\bar{h}^{\prime}(\cdot,X^{t,x;\alpha,\beta}))$. According to \cite{HH1},
Theorem 3.7, it has a unique solution in $\mathcal{S}^{2}\times\mathcal{L}%
^{2,d}\times\mathcal{M}\times\mathcal{M}$, denoted%
\[
\left(  y^{p;t,x;\alpha,\beta},z^{p;t,x;\alpha,\beta},k^{+,p;t,x;\alpha,\beta
},k^{-,p;t,x;\alpha,\beta}\right)  .
\]
By the comparison result, it is obvious that the sequence $\left(
y^{p;t,x;\alpha,\beta}\right)  _{p\in\mathbb{N}^{\ast}}$ is non-increasing.

\begin{proposition}
For every $\left(  t,x\right)  \in\left[  0,T\right]  $, $\left(  \alpha
,\beta\right)  \in\mathcal{A}_{t}\times\mathcal{B}_{t}$, $\lim_{p\rightarrow
\infty}y_{s}^{p;t,x;\alpha,\beta}=y_{s}^{t,x;\alpha,\beta}$, for all
$s\in\left[  t,T\right]  $, $\mathbb{P}$-a.s.
\end{proposition}

The proof of this result follows the same five steps of the proof of the
Theorem 3.1 in \cite{HH2} since the different approximation sequence is
constituted also by generators which are bounded and continuous. So we skip
the proof.\medskip

Let us consider, for any $p\in\mathbb{N}^{\ast}$, and any given control
processes $\alpha\left(  \cdot\right)  \in\mathcal{A}_{t}$, $\beta\left(
\cdot\right)  \in\mathcal{B}_{t}$, the associated cost functional%
\[
j^{p}\left(  t,x;\alpha,\beta\right)  :=y_{t}^{p;t,x;\alpha,\beta}%
,\quad\left(  t,x\right)  \in\left[  0,T\right]  \times\mathbb{R}^{n}%
\]
and define the lower value function\textit{ }of the approximative stochastic
differential game
\[
w^{p}\left(  t,x\right)  :=\operatorname*{essinf}_{S_{2}\in\mathbb{B}_{t}%
}\operatorname*{esssup}_{\alpha\in\mathcal{A}_{t}}j^{p}\left(  t,x;\alpha
,S_{2}\left(  \alpha\right)  \right)  .
\]

It is known, from \cite{BL2}, Proposition 3.1, that $w^{p}$ is deterministic
and is the unique viscosity solution of the equation%
\begin{equation}
\left\{
\begin{array}
[c]{l}%
\min\left\{  u\left(  t,x\right)  -\bar{h}\left(  t,x\right)  ,\max\left\{
-\dfrac{\partial u}{\partial t}\left(  t,x\right)  -H^{p}\left(
t,x,u,Du,D^{2}u\right)  ,u\left(  t,x\right)  -\bar{h}^{\prime}\left(
t,x\right)  \right\}  \right\}  =0;\\
u\left(  T,x\right)  =\bar{g}\left(  x\right)  ,
\end{array}
\right.  \label{Isaac_p}%
\end{equation}
where%
\begin{multline*}
H^{p}\left(  t,x,u,q,X\right) \\
:=\sup_{\alpha\in A}\inf_{\beta\in B}\left\{  \dfrac{1}{2}\operatorname*{Tr}%
\left(  \sigma\sigma^{T}\left(  t,x,\alpha,\beta\right)  X\right)
+\left\langle b\left(  t,x,\alpha,\beta\right)  ,q\right\rangle +f^{p}\left(
t,x,u,q\sigma\left(  t,x,\alpha,\beta\right)  ,\alpha,\beta\right)  \right\}
.
\end{multline*}

On the other hand, since the processes $\left(  y^{p;t,x;\alpha,\beta}\right)
_{p\in\mathbb{N}^{\ast}}$ form a non-increasing sequence, the sequence
$\left(  w^{p}\right)  _{p\in\mathbb{N}^{\ast}}$ is also non-increasing. By
boundedness, it has a limit, $w^{0}:\left[  0,T\right]  \times\mathbb{R}%
^{n}\longrightarrow\mathbb{R}$.

\begin{proposition}
\label{Prop subsol}The function $w^{0}$ is a viscosity subsolution of the
equation%
\begin{equation}
\left\{
\begin{array}
[c]{l}%
\min\left\{  u\left(  t,x\right)  -\bar{h}\left(  t,x\right)  ,\max\left\{
-\dfrac{\partial u}{\partial t}\left(  t,x\right)  -\bar{H}\left(
t,x,u,Du,D^{2}u\right)  ,u\left(  t,x\right)  -\bar{h}^{\prime}\left(
t,x\right)  \right\}  \right\}  =0;\\
u\left(  T,x\right)  =\bar{g}\left(  x\right)  ,
\end{array}
\right.  \label{Isaac_transf}%
\end{equation}
where%
\begin{multline*}
\bar{H}\left(  t,x,u,q,X\right) \\
:=\sup_{\alpha\in A}\inf_{\beta\in B}\left\{  \dfrac{1}{2}\operatorname*{Tr}%
\left(  \sigma\sigma^{T}\left(  t,x,\alpha,\beta\right)  X\right)
+\left\langle b\left(  t,x,\alpha,\beta\right)  ,q\right\rangle +f\left(
t,x,u,q\sigma\left(  t,x,\alpha,\beta\right)  ,\alpha,\beta\right)  \right\}
.
\end{multline*}

\end{proposition}

\begin{proof}
It is clear that $w^{0}$ is an upper semicontinuous function satisfying
$w^{0}\left(  T,\cdot\right)  \leq g\left(  \cdot\right)  $. Let now suppose
that $\varphi\in C^{1,2}\left(  [0,T)\times\mathbb{R}^{n}\right)  $ and that
$\left(  t,x\right)  \in\lbrack0,T)\times\mathbb{R}^{n}$ is a strict local
maximum point for $w^{0}-\varphi$. Then there exists a sequence $\left(
t_{p},x_{p}\right)  $ in $[0,T)\times\mathbb{R}^{n}$, converging to $\left(
t,x\right)  $, such that $w^{p}-\varphi$ has a local maximum point in $\left(
t_{p},x_{p}\right)  $ for all $p\in\mathbb{N}^{\ast}$ and $\lim_{p\rightarrow
\infty}w^{p}\left(  t_{p},x_{p}\right)  =w^{0}\left(  t,x\right)
$.$\smallskip$

Since $w^{p}$ is a viscosity solution for equation (\ref{Isaac_p}), it follows
that, for all $p\in\mathbb{N}^{\ast}$, in $\left(  t_{p},x_{p}\right)  $,%
\[
\min\left\{  w^{p}-\bar{h},\max\left\{  -\dfrac{\partial\varphi}{\partial
t}-H^{p}\left(  t_{p},x_{p},w^{p},D\varphi,D^{2}\varphi\right)  ,w^{p}-\bar
{h}^{\prime}\right\}  \right\}  =0.
\]
Because $\bar{h}\left(  t,x\right)  \leq w^{0}\left(  t,x\right)  \leq\bar
{h}^{\prime}\left(  t,x\right)  $, we analyse just two cases:

\begin{enumerate}
\item[(i)] $w^{0}\left(  t,x\right)  =\bar{h}\left(  t,x\right)  $, hence
equation (\ref{Isaac_transf}) is trivially satisfied;

\item[(ii)] $w^{0}\left(  t,x\right)  >\bar{h}\left(  t,x\right)  $, which
implies $w^{p}\left(  t_{p},x_{p}\right)  >\bar{h}\left(  t_{p},x_{p}\right)
$ for sufficiently large $p$, and so%
\[
-\dfrac{\partial\varphi}{\partial t}\left(  t_{p},x_{p}\right)  -H^{p}\left(
t_{p},x_{p},w^{p}\left(  t_{p},x_{p}\right)  ,D\varphi\left(  t_{p}%
,x_{p}\right)  ,D^{2}\varphi\left(  t_{p},x_{p}\right)  \right)  \leq0.
\]
Since%
\begin{multline*}
\dfrac{1}{2}\operatorname*{Tr}\left(  \sigma\sigma^{T}\left(  t_{p}%
,x_{p},\alpha,\beta\right)  D^{2}\varphi\left(  t_{p},x_{p}\right)  \right) \\
+\left\langle b\left(  t_{p},x_{p},\alpha,\beta\right)  ,D\varphi\left(
t_{p},x_{p}\right)  \right\rangle +f^{p}\left(  t_{p},x_{p},w^{p}\left(
t_{p},x_{p}\right)  ,D\varphi\left(  t_{p},x_{p}\right)  \sigma\left(
t_{p},x_{p},\alpha,\beta\right)  ,\alpha,\beta\right)
\end{multline*}
converges uniformly (with respect to $\alpha$ and $\beta$) to%
\begin{multline*}
\dfrac{1}{2}\operatorname*{Tr}\left(  \sigma\sigma^{T}\left(  t,x,\alpha
,\beta\right)  D^{2}\varphi\left(  t,x\right)  \right) \\
+\left\langle b\left(  t,x,\alpha,\beta\right)  ,D\varphi\left(  t,x\right)
\right\rangle +\tilde{f}\left(  t,x,w^{0}\left(  t,x\right)  ,D\varphi\left(
t,x\right)  \sigma\left(  t,x,\alpha,\beta\right)  ,\alpha,\beta\right)  ,
\end{multline*}
and $w^{0}\left(  t,x\right)  \leq M$, it follows that%
\[
-\dfrac{\partial\varphi}{\partial t}\left(  t,x\right)  -\bar{H}\left(
t,x,w^{0}\left(  t,x\right)  ,D\varphi\left(  t,x\right)  ,D^{2}\varphi\left(
t,x\right)  \right)  \leq0.
\]

\end{enumerate}

This finishes our proof.\medskip
\end{proof}

One can repeat the above schema, but with lower approximation, \emph{i.e.} we
can construct (in the same manner), an increasing sequence of Lipschitz,
bounded functions $\left(  f_{p}\right)  $ converging to $\tilde{f}$. By
denoting%
\[
(y_{p}^{t,x;\alpha,\beta},z_{p}^{t,x;\alpha,\beta},k_{p}^{+,t,x;\alpha,\beta
},k_{p}^{-,t,x;\alpha,\beta})
\]
the solution of the RBSDE associated with $(f_{p},\bar{g}(X_{T}^{t,x;\alpha
,\beta}),\bar{h}(\cdot,X^{t,x;\alpha,\beta}),\bar{h}^{\prime}(\cdot
,X^{t,x;\alpha,\beta}))$, one can show that $\left(  y_{p}^{t,x;\alpha,\beta
}\right)  $ is an increasing sequence of processes, converging to a minimal
solution of the RBSDE associated with $(f,\bar{g}(X_{T}^{t,x;\alpha,\beta
}),\bar{h}(\cdot,X^{t,x;\alpha,\beta}),\bar{h}^{\prime}(\cdot,X^{t,x;\alpha
,\beta}))$. Let us denote%

\[
j_{p}\left(  t,x;\alpha,\beta\right)  :=y_{p,t}^{t,x;\alpha,\beta}%
,\quad\left(  t,x\right)  \in\left[  0,T\right]  \times\mathbb{R}^{n}%
\]
and
\[
w_{p}\left(  t,x\right)  :=\operatorname*{essinf}_{S_{2}\in\mathbb{B}_{t}%
}\operatorname*{esssup}_{\alpha\in\mathcal{A}_{t}}j_{p}\left(  t,x;\alpha
,S_{2}\left(  \alpha\right)  \right)  .
\]
Then, setting $w_{0}:=\lim w_{p}$, we have the analogous of Proposition
\ref{Prop subsol}, whose proof is essentially the same.

\begin{proposition}
The function $w_{0}$ is a viscosity supersolution of the equation
(\ref{Isaac_transf}).
\end{proposition}

Let us now define $\mathcal{W}^{0}:=\ln\dfrac{w^{0}}{2C}$ and $\mathcal{W}%
_{0}:=\ln\dfrac{w_{0}}{2C}$. It is straightforward to show that $\mathcal{W}%
^{0}$ and $\mathcal{W}_{0}$ are viscosity subsolution, respectively
supersolution, of equation (\ref{Isaac1}). By the comparison result (see
\cite{HH2}, Remark 3.3), $\mathbb{P}$-a.s.,%
\[
\mathcal{W}^{0}\left(  t,x\right)  \geq\mathcal{W}\left(  t,x\right)
\geq\mathcal{W}_{0}\left(  t,x\right)  ,\ \forall\left(  t,x\right)
\in\left[  0,T\right]  \times\mathbb{R}^{n},
\]
since, for all $p\in\mathbb{N}^{\ast}$, $\left(  t,x\right)  \in\left[
0,T\right]  \times\mathbb{R}^{n}$ and $\left(  \alpha\left(  \cdot\right)
,\beta\left(  \cdot\right)  \right)  \in\mathcal{A}_{t}\times\mathcal{B}_{t}$,
$\mathbb{P}$-a.s.,%
\[
y_{s}^{p;t,x;\alpha,\beta}\geq y_{s}^{t,x;\alpha,\beta}\geq y_{p,s}%
^{t,x;\alpha,\beta},\ \forall s\in\left[  t,T\right]  .
\]
On the other hand, by Theorem \ref{uniqueness}, we have that $\mathcal{W}%
^{0}\leq\mathcal{W}_{0}$. These inequalities imply that $\mathcal{W}$ is
deterministic and is a viscosity solution of equation (\ref{Isaac1}). The
existence result is thus proved.

\begin{remark}
As one can see from the proof of the existence result, that $\mathcal{W}$ can
be defined via any solution of the BSDE associated with $\left(
F,g,h,h^{\prime}\right)  $, not necessarily the maximal one. This could be
used to deduce the uniqueness for this equation.
\end{remark}

\section{Applications}

As we said in the introduction, one application comes from financial markets.
Indeed, it is possible to use our framework for the study of American game
options and of the Ramsey's model. We now focus on the link between the
payoffs $\Gamma$ defined in (\ref{payoff}) and the solutions of the reflected
BSDEs associated with (\ref{mainRBSDE}).

\begin{proposition}
Assume that $F(t,x,y,z,\alpha,\beta)=\varphi(t,x,\alpha,\beta)+\frac{1}%
{2}|z|^{2}$, where $\varphi$ is a bounded measurable function. For any
stopping times $\tau$ and $\sigma$ let us consider the following standard
BSDE:%
\begin{equation}
\left\{
\begin{array}
[c]{l}%
(Y^{(t,x);(\alpha,\sigma),(\beta,\tau)},Z^{(t,x);(\alpha,\sigma),(\beta,\tau
)})\in\mathcal{S}^{2}\times\mathcal{H}^{2,d}\medskip\\
Y_{s}^{(t,x);(\alpha,\sigma),(\beta,\tau)}=h(X_{\sigma}^{t,x;\alpha,\beta
})1_{[\sigma\leq\tau<T]}+h^{\prime}(X_{\tau}^{t,x;\alpha,\beta})1_{[\tau
<\sigma]}+g(X_{T}^{t,x;\alpha,\beta})1_{[\sigma=\tau=T]}\medskip\\
\qquad\qquad\qquad\qquad\qquad+%
{\displaystyle\int_{s\wedge\tau\wedge\sigma}^{T\wedge\tau\wedge\sigma}}
\left\{  \varphi(r,X_{r}^{t,x;\alpha,\beta},\alpha_{r},\beta_{r})+\dfrac{1}%
{2}|Z_{r}^{(t,x);(\alpha,\sigma),(\beta,\tau)}|^{2}\right\}  dr\medskip\\
\qquad\qquad\qquad\qquad\qquad\qquad\qquad\qquad-%
{\displaystyle\int_{s\wedge\tau\wedge\sigma}^{T\wedge\tau\wedge\sigma}}
Z_{r}^{(t,x);(\alpha,\sigma),(\beta,\tau)}dW_{r},\,\forall s\leq T.
\end{array}
\right.  \label{application eq}%
\end{equation}
Then
\[
\Gamma(\alpha,\sigma;\beta,\tau)=\exp\{Y_{0}^{(0,x);(\alpha,\sigma
),(\beta,\tau)}\},\smallskip
\]%
\[
\operatorname*{essinf}_{S_{2}\in\mathbb{B}_{t}}\operatorname*{esssup}%
_{\alpha\in\mathcal{A}_{t}}Y_{t}^{t,x;\alpha,S_{2}(\alpha)}%
=\operatorname*{essinf}_{S_{2}\in\mathbb{B}_{t}}\operatorname*{esssup}%
_{\alpha\in\mathcal{A}_{t}}\operatorname*{essinf}_{\tau\in\mathcal{T}_{t}%
}\operatorname*{esssup}_{\sigma\in\mathcal{T}_{t}}Y_{t}^{(t,x);(\alpha
,\sigma),(S_{2}\left(  \alpha\right)  ,\tau)}%
\]
and
\[
\operatorname*{esssup}_{S_{1}\in\mathbb{A}_{t}}\operatorname*{essinf}%
_{\beta\in\mathcal{B}_{t}}Y_{t}^{t,x;S_{1}(\beta),\beta}%
=\operatorname*{esssup}_{S_{1}\in\mathbb{A}_{t}}\operatorname*{essinf}%
_{\beta\in\mathcal{B}_{t}}\operatorname*{esssup}_{\sigma\in\mathcal{T}_{t}%
}\operatorname*{essinf}_{\tau\in\mathcal{T}_{t}}Y_{t}^{(t,x);(S_{1}\left(
\beta\right)  ,\sigma),(\beta,\tau)},
\]
where by $\mathcal{T}_{t}$ we denoted the set of stopping times $\tau$ such
that $t\leq\tau\leq T$.
\end{proposition}

\begin{proof}
For the sake of simplicity we denote $Y^{(t,x);(\alpha,\sigma),(\beta,\tau)}$
by $Y$. First note that since the functions $h$, $h^{\prime}$, $g$ and
$\varphi$ are bounded then through the result by Kobylanski \cite{K}, Theorem
2.3, the solution of (\ref{application eq}) exists and is unique. Now for
$s\leq T$ let us set%
\[
\tilde{Y}_{s}=\exp\left\{  Y_{s}+\int_{0}^{s\wedge\tau\wedge\sigma}%
\varphi(r,X_{r}^{t,x;\alpha,\beta},\alpha_{r},\beta_{r})dr\right\}
\]
Using now It\^{o}'s formula to obtain that:%
\[
\left\{
\begin{array}
[c]{l}%
d\tilde{Y}_{s}=\tilde{Z}_{s}dW_{s},\,\,s\leq T\mbox{ and }\smallskip\\
\tilde{Y}_{T}=\exp\left\{  h(X_{\sigma}^{t,x;\alpha,\beta})1_{[\sigma\leq
\tau<T]}+h^{\prime}(\tau,X_{\tau}^{t,x;\alpha,\beta})1_{[\tau<\sigma]}%
+g(X_{T}^{t,x;\alpha,\beta})1_{[\sigma=\tau=T]}\right. \\
\quad\quad\quad\quad\quad\quad\quad\quad\quad\quad\quad\quad\quad\quad
\quad\quad\quad\quad\quad\left.  +%
{\displaystyle\int_{0}^{T\wedge\tau\wedge\sigma}}
\varphi(r,X_{r}^{t,x;\alpha,\beta},\alpha_{r},\beta_{r})dr\right\}  .
\end{array}
\right.
\]
It implies that $\mathbb{E}[\tilde{Y}_{0}]=\mathbb{E}[\tilde{Y}_{T}]$. As
$\tilde{Y}_{0}$ is deterministic since it is $\mathcal{F}_{0}$-measurable then
$\tilde{Y}_{0}=\mathbb{E}[\tilde{Y}_{0}]=\mathbb{E}[\tilde{Y}_{T}]$, i.e.
$\exp\{Y_{0}\}=\Gamma(\alpha,\sigma;\beta,\tau)$.\medskip

Let us now deal with the second relations. First note (see for example
\cite{elkal},
Proposition 2.3%
) that the characterization of a solution for a BSDE with two reflecting
barriers implies that for any $\alpha$ and $\beta$,%
\begin{equation}
Y_{t}^{t,x;\alpha,\beta}=\operatorname*{essinf}_{\tau\geq t}%
\operatorname*{esssup}_{\sigma\geq t}Y_{t}^{(t,x);(\alpha,\sigma),(\beta
,\tau)}=\operatorname*{esssup}_{\sigma\geq t}\operatorname*{essinf}_{\tau\geq
t}Y_{t}^{(t,x);(\alpha,\sigma),(\beta,\tau)}\label{caracterization}%
\end{equation}
since%
\[
\exp(Y_{t}^{t,x;\alpha,\beta})=\operatorname*{essinf}_{\tau\geq t}%
\operatorname*{esssup}_{\sigma\geq t}\exp\{Y_{t}^{(t,x);(\alpha,\sigma
),(\beta,\tau)}\}=\operatorname*{esssup}_{\sigma\geq t}\operatorname*{essinf}%
_{\tau\geq t}\exp\{Y_{t}^{(t,x);(\alpha,\sigma),(\beta,\tau)}\}.
\]
Next, for any $\alpha\in\mathcal{A}_{t}$ and $S_{2}\in\mathbb{B}_{t}$, the
formula
(\ref{caracterization}) implies clearly that the first equality holds. The
second one is treated in the same manner.%

\end{proof}

\begin{remark}
The following relation holds: $\forall s\leq T$, $\forall\tau,\sigma
\in\mathcal{T}_{t}$,
\[%
\begin{array}
[c]{l}%
Y_{s}^{(t,x);(\alpha,\sigma),(\beta,\tau)}=\ln\left(  \mathbb{E}\left[
\exp\{h(X_{\sigma}^{t,x;\alpha,\beta})1_{[\sigma\leq\tau<T]}+h^{\prime
}(X_{\tau}^{t,x;\alpha,\beta})1_{[\tau<\sigma]}+g(X_{T}^{t,x;\alpha,\beta
})1_{[\sigma=\tau=T]}\right.  \right.  \medskip\\
\qquad\qquad\qquad\qquad\qquad\qquad\qquad\qquad\qquad\qquad\left.  \left.
\left.  +\int_{s}^{T\wedge\tau\wedge\sigma}\varphi(r,X_{r}^{t,x;\alpha,\beta
},\alpha_{r},\beta_{r})dr\}\right\vert \mathcal{F}_{s}\right]  \right)
.\medskip
\end{array}
\]

\end{remark}

\begin{remark}
Of course, a more interesting and difficult issue is when the upper and the
lower values of the mixed zero-sum two-players stochastic differential game
are equal respectively to
\[
\inf_{\left(  \beta,\tau\right)  }\sup_{\left(  \alpha,\sigma\right)  }%
\Gamma(\alpha,\sigma;\beta,\tau)\quad\mbox{ and }\quad\sup_{\left(
\alpha,\sigma\right)  }\inf_{\left(  \beta,\tau\right)  }\Gamma(\alpha
,\sigma;\beta,\tau).
\]
Unfortunately, this is still an open problem.
\end{remark}


\bigskip

\begin{thebibliography}{99}                                                                                               %


\bibitem {ab}Amilon H.; Bermin H.P. - \textit{Welfare effects of controlling
labor supply: an application of the stochastic Ramsey model}, Journal of
Economic Dynamics \& Control 28, pp.331-348 (2003).

\bibitem {BHM}Bahlali, S.; Hamad\`{e}ne, S.; Mezerdi, B. - \textit{Backward
stochastic differential equations with two reflecting barriers and quadratic
growth coefficient}, Stochastic Processes and Their Applications 115, no.7,
pp.1107-1129 (2005).

\bibitem {bk}Baurdoux, C.; Kyprianou, A. - \textit{Further Calculations for
Israeli Options,} Stochastics.and Stochastics Reports., 76, pp. 549-569 (2004).

\bibitem {bensnagai}Bensoussan, A.; Nagai, H. - \textit{Min-max
characterization of a small noise limit on risk-sensitive control}, SIAM J.
Control Optim. 35 (4), pp. 1093--1115 (1997).

\bibitem {bensnagai2}Bensoussan, A.; Frehse, J.; Nagai, H. - \textit{Some
results on risk-sensitive with full observation}, J. Appl. Math. Optim. 37 (1998).

\bibitem {bcr}Buckdahn, R.; Cardialaguet, P.; Rainer, C. - \textit{Nash
equilibrium payoffs for nonzero-sum stochastic differential games}, SIAM J.
Cont. Opt. 43, No.2, 624-642. 35 (2004).

\bibitem {BL1}Buckdahn, R.; Li J. - \textit{Stochastic Differential Games with
Reflection and Related Obstacle Problems for Isaacs Equations},
http:/arXiv:0707.1133v2 [math.PR] (25 Jul 2007).

\bibitem {BL2}Buckdahn, R.; Li J. - \textit{Probabilistic Interpretation for
Systems of Isaacs Equations with Two Reflecting Barriers}, available at
http:/arXiv:0804.0311v1 [math.OC] (2 Apr 2008).

\bibitem {CIL}Crandall, M.G.; Ishii, H.; Lions, P.L. - \textit{User's Guide to
Viscosity Solutions of Second Order Partial Differential Equations}, Bulletin
of The American Mathematical Society, Volume 27, Number 1, July 1992, pp 1-67.

\bibitem {CK}Cvitanic, J.; Karatzas, I. - \textit{Backward stochastic
differential equations with reflection and Dynkin games}, The Annals of
Probability 24, no. 4, pp. 2024-2056 (1996).

\bibitem {dupuis}Dupuis, P.; McEneaney W.M. - {\textit{Risk-sensitive and
robust escape criteria}, SIAM J. Control Optim. 35 (6), 2021--2049 (1996). }

\bibitem {snic}N. El-Karoui; S. Hamad\`{e}ne - \textit{BSDEs and
risk-sensitive control, zero-sum and nonzerosum game problems of stochastic
functional differential equations}, Stochastic Process Appl., 107, pp.
145--169 (2003).

\bibitem {elkal}N. El-Karoui; C. Kapoudjian; E. Pardoux; S. Peng; M. C. Quenez
\textit{ - Reflected solutions of backward SDEs and related obstacle problems
for PDEs}, Ann. Probab., 25, pp. 702--737 (1997).

\bibitem {FS}Fleming, W.H.; Souganidis, P.E. - \textit{On the existence of
value functions of two-player, zero-sum stochastic differential games},
Indiana Univ. Math. J. 38, No. 2, pp. 293-314 (1989).

\bibitem {H}Hamad\`{e}ne, S - \textit{Mixed zero-sum differential game and
American game options}, SIAM JCO, Vol. 45 (2), pp. 496-518 (2006)

\bibitem {HH1}Hamad\`{e}ne, S.; Hassani M. - \textit{BSDEs with two reflecting
barriers: the general result}, Probab. Theory Relat. Fields 132, pp. 237-264 (2005).

\bibitem {HH2}Hamad\`{e}ne, S.; Hdhiri I. - \textit{Backward stochastic
differential equations with two distinct reflecting barriers and quadratic
growth generator}, Journal of Applied Mathematics and Stochastic Analysis,
Vol. 2006, Article ID 95818, 28 pages.

\bibitem {hl1}Hamad\`{e}ne, S.; Lepeltier, J.P. - \textit{Zero-sum stochastic
differential games and backward equations,} Systems and Control Letters. 24,
259-263 (1995).

\bibitem {hl2}Hamad\`{e}ne, S.; Lepeltier,.J.P. - \textit{Backward equations,
stochastic control and zero-sum stochastic differential games}, Stochastics
and Stochastic Reports, vol.54, pp.221-231 (1995).

\bibitem {hl3}Hamad\`{e}ne, S.; Lepeltier, J.P. - \textit{Reflected BSDEs and
mixed game problems}, Stochastic processes and their applications, 85, 177-188 (2000).

\bibitem {hlp}Hamad\`{e}ne, S.; Lepeltier, J.P.; Peng, S. - \textit{BSDEs with
continuous coefficients and stochastic differential games,} El Karoui, N. and
Mazliak, L. (Eds.), Backward stochastic differential equations. Harlow:
Longman. Pitman Res. Notes Math. Ser. 364, 115-128 (1997).

\bibitem {hlw}Hamad\`{e}ne, S.; Lepeltier, J.P.; Wu, Z. - \textit{Infinite
horizon Reflected BSDEs and applications in mixed control and game problems,}
Probability and mathematical statistics, 19, 211-234 (1999).

\bibitem {kk}Kallsen, J.; Kuhn, C. - \textit{Pricing Derivatives of American
and game type in Incomplete Markets}, Finance Statistics (2003).

\bibitem {yk}Y. Kifer - \textit{Game options}, Finance Stoch., 4, pp. 443--463 (2000).

\bibitem {K}Kobylanski, M. - \textit{Backward Stochastic Differential
Equations and Partial Differential Equations with Quadratic Growth}, The
Annals of Probability, Vol. 28, No. 2, pp. 558-602 (2000).

\bibitem {kk2}Kyprianou, A.; Kuhn, C. - \textit{Israeli Options as composite
Exotic Options} (2003).

\bibitem {KS}Karatzas, I.; Shreve, S.E. - \textit{Brownian motion and
Stochastic Calculus}, Springer-Verlag, N.Y. (1991).

\bibitem {nh}Nagai, H. - \textit{Bellmann equations of risk-sensitive
control}, SIAM J. Control Optim. 34 (1), pp. 74--101 (1996).
\end{thebibliography}
\end{document}